\documentclass[12pt]{article}
\usepackage{a4wide}
\usepackage{amsmath}
\usepackage{amsfonts}
\usepackage{amsthm}
\usepackage{amssymb}
\theoremstyle{plain}
\newtheorem{theorem}{Theorem}[section] 
\newtheorem*{theorem*}{Theorem}  
\newtheorem{proposition}[theorem]{Proposition}

\newtheorem{definition}[theorem]{Definition}

\theoremstyle{remark}

\newtheorem*{example*}{Example}

\numberwithin{equation}{section}

\newcommand{\Z}{\mathbb{Z}}
\newcommand{\Q}{\mathbb{Q}}
\newcommand{\R}{\mathbb{R}}
\newcommand{\C}{\mathbb{C}}

\renewcommand{\H}{\mathbb{H}}
\newcommand{\leg}[2]{\left( \frac{#1}{#2} \right)}
\newcommand{\zxz}[4]{\begin{pmatrix} #1 & #2 \\ #3 & #4 \end{pmatrix}}
\newcommand{\abcd}{\zxz{a}{b}{c}{d}}

\newcommand{\kzxz}[4]{\left(\begin{smallmatrix} #1 & #2 \\ #3 & #4\end{smallmatrix}\right) }
\newcommand{\kabcd}{\kzxz{a}{b}{c}{d}}

\newcommand{\calC}{{\cal C}}

\newcommand{\calE}{{\cal E}}

\newcommand{\calG}{{\cal G}}
\newcommand{\calH}{{\cal H}}

\newcommand{\calK}{{\cal K}}
\newcommand{\calL}{{\cal L}}

\newcommand{\calO}{{\cal O}}

\newcommand{\calW}{{\cal W}}

\newcommand{\fraka}{\mathfrak a}
\newcommand{\frakd}{\mathfrak d}
\newcommand{\frake}{\mathfrak e}

\newcommand{\frakg}{\mathfrak g}
\newcommand{\frakp}{\mathfrak p}
\newcommand{\frakq}{\mathfrak q}
\newcommand{\frakk}{\mathfrak k}
\newcommand{\frakh}{\mathfrak h}

\newcommand{\bs}{\backslash}
\newcommand{\OK}{\calO}
\renewcommand{\div}{\operatorname{div}}
\newcommand{\norm}{\operatorname{N}}
\newcommand{\vol}{\operatorname{vol}}

\newcommand{\sgn}{\operatorname{sgn}}

\newcommand{\dv}{\operatorname{div}}

\newcommand{\Sl}{\operatorname{SL}}
\newcommand{\Symp}{\operatorname{Sp}}
\newcommand{\Mp}{\operatorname{Mp}}
\newcommand{\Orth}{\operatorname{O}}

\newcommand{\Gr}{\operatorname{Gr}}

\begin{document}

\title{Integrals of automorphic Green's functions associated to Heegner divisors}
\date{April 10, 2002\thanks{revision December 10, 2002}}
\author{Jan Hendrik Bruinier\thanks{Supported by a Heisenberg-Stipendium of the DFG.}\; and Ulf K\"uhn}
\maketitle

\noindent

\subsection*{Abstract}
In the present paper we find explicit formulas for the degrees 
of Heegner divisors on arithmetic quotients of the orthogonal 
group $\Orth(2,p)$ and for the integrals of certain automorphic 
Green's functions associated with Heegner divisors. The latter 
quantities are important in the study of the arithmetic degrees 
of Heegner divisors in the context of Arakelov geometry.
In particular, we obtain a different proof 
and a generalization of 
results of Kudla relating 
these quantities to the Fourier coefficients of 
certain non-holomorphic Eisenstein series of weight $1+p/2$ for the
 metaplectic group $\Mp_2(\Z)$.

\section{Introduction}

Integrals of automorphic Green's functions contribute to the
calculation of arithmetic degrees of regular models of Shimura
varieties over Dedekind rings 
and Faltings heights of its subvarieties. 
These integrals are expected to be related to logarithmic derivatives 
of certain 
$L$-functions.
For example, if $E_4(\tau)= 1+240 \sum_{n>0} \sigma_{3}(n) e^{2 \pi i n \tau}$, with
$\tau=x+iy$ in the complex upper half plane $\H$, is
the Eisenstein series of weight $4$ for the modular group
$\Sl_2(\Z)$, then 
\begin{align}\label{eq:werbung}
\frac{1}{L(\chi_{-3},0)} 
\int\limits_{\Sl_2(\Z)\setminus \H} \log 
\left( |E_4(\tau)|^2(4 \pi y)^4 \right)
\frac{dx\,dy}{4 \pi y^2}= 
2  \frac{\zeta'(-1)}{\zeta(-1) } + 1 -
  \frac{L'(\chi_{-3},0)}{L(\chi_{-3},0)} - \frac{1}{2} \log(3) ,
\end{align}
where $L(\chi_{-3},s)$ and $\zeta(s)$ denote the Dirichlet
$L$-function for the character $\chi_{-3}$ and the Riemann zeta
function, respectively. Notice that $L(\chi_{-3},0)=1/3$.
The core of such formulas is the Kronecker limit formula.  

In the present paper we generalize the above identity to 
modular forms on the orthogonal group $\Orth(2,p)$ that vanish on 
Heegner divisors.
More precisely, we derive an explicit formula expressing the integrals of 
certain Green's function associated to Heegner divisors on $\Orth(2,p)$ in 
terms of logarithmic derivatives of Dirichlet $L$-functions, the Riemann zeta 
function, and generalized divisor sums.  
As a corollary, we obtain a generalization of \eqref{eq:werbung} to Borcherds 
products on $\Orth(2,p)$ \cite{Bo1,Bo2}. 
Observe that $E_4$ in the above example can be viewed as a Borcherds product. Its 
divisor on $\Sl_2(\Z)\bs\H$ is given by the Heegner point $e^{\pi i/3}$.  
Our approach relies on the construction of Borcherds products using automorphic
 Green's functions introduced in \cite{Br1}, which can be viewed as a generalization
 of the Kronecker limit formula; combined with an integral identity for 
such Green's functions as in \cite{OT}. 
In  \cite{Ku2} S.~Kudla calculated the integral of the logarithm of the Petersson norm of Borcherds products using 
the Siegel-Weil formula together with a Stoke's type argument. 
We discuss the relation of his results to ours below.

\bigskip

We now describe the content of this paper in more detail.
Let $(V,q)$ be a real quadratic space of signature $(2,p)$ and put $\kappa=1+p/2$. 
We assume that either $p> 2$, or that $p=2$ and the dimension of a maximal isotropic subspace of $L\otimes_\Z\Q$ equals $1$. 
We write $\calK'$ for the irreducible Hermitean symmetric space of dimension $p$ associated with the real orthogonal group of $(V,q)$. 
Let $L\subset V$ be an even lattice, and $L'$ its dual. We denote by $\Gamma(L)$ the discriminant kernel of the orthogonal group $\Orth(L)$ of $L$, that is the kernel of the natural homomorphism from $\Orth(L)$ to $\Orth(L'/L)$.
By the theory of Baily-Borel the arithmetic quotient 
\[
X_L=\Gamma(L)\bs \calK'
\]
is a quasi-projective algebraic variety. 
Let $\Omega$ be the K\"ahler form on $X_L$ given by the first Chern class of the line bundle of modular forms of weight $1$ on $X_L$.
We define the volume of $X_L$ by $\vol(X_L)=\int_{X_L}\Omega^p$ and the degree of a divisor $D$ on $X_L$  
by
\[
\deg(D)=\int_{D}\Omega^{p-1}.
\]

Recall that for any $\beta\in L'/L$ and any negative $m\in \Z+q(\beta)$
there is a certain special divisor $H(\beta,m)$ on $X_L$, called the Heegner divisor of discriminant $(\beta,m)$ (see section \ref{sect3} for a precise definition). These divisors arise from embedded quotients analogous to $X_L$ of dimension $p-1$. They include Heegner points on modular curves, Hirzebruch-Zagier divisors on Hilbert modular surfaces, and Humbert surfaces on Siegel modular threefolds as special cases 
and have been studied by many people, e.g.~\cite{Bo3}, \cite{KM}, \cite{Ku:DUKE}, \cite{Oda:SIEGEL}, \cite{Br1}.

For any Heegner divisor $H(\beta,m)$ there exists an associated automorphic Green's function $\Phi_{\beta,m}(Z,s)$, where $Z\in \calK'$ and $s\in \C$ with $\Re(s)>\kappa/2$. As a function in $Z$ it is an eigenfunction  of the invariant Laplacian on $\calK'$ and has a logarithmic singularity along $H(\beta,m)$.
Such Green's functions were introduced in \cite{Br1,Br2}, and independently from a different perspective in \cite{OT}.

Using the approach of \cite{Br1} the Fourier expansions of the functions  $\Phi_{\beta,m}(Z,s)$ can be determined. Moreover, one can show that they have a meromorphic continuation in $s$ to a neighborhood of $\kappa/2$ with a simple pole at $s=\kappa/2$. It turns out that their singularities at $s=\kappa/2$ are dictated by the coefficients of a certain non-holomorphic Eisenstein series of weight $\kappa/2$ for $\Mp_2(\Z)$. 
We describe this in somewhat more detail. To simplify the exposition we temporarily assume that $p$ is even (therefore $\kappa\in \Z$) and $p>2$. In the general case, treated in the body of this paper, one has to replace $\Sl_2(\Z)$ by the metaplectic group $\Mp_2(\Z)$. 

Let $\rho_L$ be the Weil representation of $\Sl_2(\Z)$ on the group algebra $\C[L'/L]$ of the (finite) discriminant group of $L$ as in \cite{Bo2}, \cite{Br1}, and denote by $(\frake_\gamma)_{\gamma\in L'/L}$ the standard basis of $\C[L'/L]$. We define the Eisenstein series $E_0(\tau,s)$ by
\[
E_0(\tau,s)= \sum_{M\in \Gamma_\infty\bs \Sl_2(\Z)} 
\frac{\Im(M\tau)^s}{j(M,\tau)^\kappa} \big(\rho_L^t(M)\frake_0\big)  ,  
\]
where $\Gamma_\infty=\{\kzxz{1}{n}{0}{1};\; n\in \Z\}$ and $j(\kabcd,\tau)=c\tau+d$.
It is a $\C[L'/L]$-valued non-holomorphic modular form of weight $\kappa$ with respect to $\Sl_2(\Z)$ and the dual of $\rho_L$. It has a Fourier expansion of the form
\begin{align*}
E_0(\tau,s) &= \sum_{\gamma\in L'/L}\sum_{n\in \Z-q(\gamma)} c_0 (\gamma,n,s,y)  e(n x) \frake_\gamma,
\end{align*}
where  $e(x)=e^{2\pi i x}$ as usual.
The Fourier coefficients with non-zero index decompose into the product
\[
c_0(\gamma,n,s,y)=C(\gamma,n,s)\calW_s(4\pi ny)
\]
of a coefficient $C(\gamma,n,s)$, which is independent of $y$, and a part given by a Whittaker function $\calW_s(y)$ (see \eqref{eq:W}).
In Theorem \ref{thm:eismain} we derive a closed formula for these coefficients following the argument of \cite{BK}, which relies on Shintani's formula \cite{Sh} for the coefficients of the Weil representation, and a result of Siegel on representation numbers of quadratic forms modulo prime powers.
This theorem is vital for all explicit computations in the present paper.

In Propositions \ref{prop:aux} and \ref{prop:pc} we show that the residue of $\Phi_{\beta,m}(Z,s)$ at $s=\kappa/2$ is equal to $-C(\beta,-m,0)$. It is a consequence of this fact and the fundamental integral formula for the Green's function $\Phi_{\beta,m}(Z,s)$ which follows from the work of Oda and Tsuzuki \cite{OT} (see Theorem \ref{thm:kron}) that the special value $E_0(\tau,0)$ is the generating series for the 
degrees of Heegner divisors. We recover the identity 
\begin{equation}\label{intro:1}
E_0(\tau,0)=2\frake_0-\frac{2}{\vol(X_L)}\sum_{\gamma\in L'/L}
\sum_{\substack{n>0}}\deg(H(\gamma,-n)) e(n\tau)\frake_\gamma,
\end{equation}
see e.g.~\cite{Ge2}, \cite{HZ}, \cite{Ku2}, \cite{Oda:SIEGEL}.
Variants and generalizations of this result were also proved by Kudla-Millson in their work on special cycles (see e.g.~\cite{Kudla:Bourbaki}).
By means of our formula for the coefficients $C(\beta,m,s)$ we 
find explicit formulas for the degrees of Heegner divisors in terms of 
special values of Dirichlet $L$-functions, the Riemann zeta function, 
and generalized divisor sums (see Proposition \ref{thm:deg}).

In \cite{Bo1, Bo2} Borcherds constructed meromorphic modular forms for the group $\Gamma(L)$ with zeros and poles on Heegner divisors as multiplicative liftings of $\C[L'/L]$-valued nearly holomorphic modular forms of weight $1-p/2$ for the group $\Mp_2(\Z)$ and the Weil representation $\rho_L$.
They have explicit infinite product expansions analogous to the Delta function.
In particular, their Fourier expansions have integral cyclotomic coefficients.

The above Green's functions and Borcherds products are related in the following way. We define the regularized Green's function $G_{\beta,m}(Z)$ to be $-1/4$ times the constant term in the Laurent expansion of $\Phi_{\beta,m}(Z,s)$ in $s$ at $s=\kappa/2$ plus some normalizing constant which essentially involves the coefficient $C(\beta,-m,s)$ and its derivative at $s=\kappa/2$.
We denote by $\|\cdot\|$ the Petersson metric on the line bundle of modular forms of weight $k$ on $X_L$ normalized as in \eqref{eq:petnorm}.
It was proved in \cite{Br2} that
if $F$ is a Borcherds product (in the sense of \cite{Bo2} Theorem 13.3) with divisor $\dv(F)=\frac{1}{2}\sum_{\beta,m} a(\beta,m) H(\beta,m)$, then 
\begin{equation}\label{intro:2}
\log\|F(Z)\|^2 = \sum_{\beta,m}a(\beta,m) G_{\beta,m}(Z).
\end{equation}
So the Green's functions $G_{\beta,m}(Z)$ can be viewed as the building blocks of Borcherds products. Observe that the individual functions $G_{\beta,m}(Z)$ are in general far from being the Petersson norm of a modular form.

\begin{theorem*}
(See Theorem \ref{thm:main}.)
The integral of $G_{\beta,m}(Z)$ is related to the logarithmic derivative of the coefficient $C(\beta,-m,s)$ of the Eisenstein series $E_0(\tau,s)$ at $s=0$ as follows:
\begin{equation}\label{intro:3}
\frac{2}{\deg H(\beta,m)}\int\limits_{X_L} G_{\beta,m}(Z) \Omega^p 
=\frac{C'(\beta,-m,0)}{C(\beta,-m,0)}+\log(4\pi)-\Gamma'(1).
\end{equation}

Moreover, if the rank $r$ of $L$ is even, then the right hand side of \eqref{intro:3} is equal to
\[
2\frac{L'(\chi_{D_0},1-\kappa)}{L(\chi_{D_0},1-\kappa)}+2\frac{\sigma_{\beta,-m}'(\kappa)}{\sigma_{\beta,-m}(\kappa)}+\log|mD_0^2|
+\sum_{j=1}^{\kappa-1}\frac{1}{j},
\]
where $D_0$ denotes the discriminant of the quadratic field $\Q(\sqrt{d})$ with  $d=(-1)^{r/2}\det(L)$. 

If $r$ is odd, then the right hand side of \eqref{intro:3} is equal to 
\[
4\frac{\zeta'(2-2\kappa)}{\zeta(2-2\kappa)}-2\frac{L'(\chi_{D_0},3/2-\kappa)}{L(\chi_{D_0},3/2-\kappa)}
+2\frac{\sigma_{\beta,-m}'(\kappa)}{\sigma_{\beta,-m}(\kappa)}+\log|4m/D_0^2|
+\sum_{j=1}^{\kappa-1/2}\frac{2}{2j-1},
\]
where $D_0$ is the discriminant of the quadratic field $\Q(\sqrt{md})$ with $d=2(-1)^{(r-1)/2} \det(L)$.

In both cases we understand by $\det(L)$ the Gram determinant of $L$, by $L(\chi_{D_0},s)$ the Dirichlet $L$-function associated with the character $\chi_{D_0}$, and by $\sigma_{\gamma,n}(s)$ the generalized divisor sum defined in \eqref{eq:defsigma}.
\end{theorem*}

The contribution coming from the logarithmic logarithmic derivative of $\sigma_{\beta,-m}(s)$ at $\kappa$ is a sum of the form $\sum_{p} \alpha_p\log(p)$ over the primes dividing $2 m \det(L)^2$, where the $\alpha_p$ are rational coefficients depending on the representation numbers of the lattice $L$ modulo powers of $p$.

In view of \eqref{intro:2}, as a Corollary, one gets explicit formulas for the integrals of the Petersson norms of Borcherds products. So our result yields a
 different proof and a generalization of the main theorems  in 
\cite{Ku2}, \cite{KuYa}. The second assertion of the theorem is
 the desired generalization of \eqref{eq:werbung} to the group $\Orth(2,p)$.
Formally we recover \eqref{eq:werbung} using the exceptional isomorphism 
between $\Orth^0(2,1)$ and $\Sl_2(\R)$ and the arithmetic subgroup defined by the lattice $L=\Z(2)\perp II_{1,1}$ of rank $3$ with Gram determinant $-2$. The logarithm of the Petersson norm of $E_4$ is equal to the Green's function for the Heegner divisor $H(\beta,-3/4)$ with $m=-3/4$ and $D_0=-3$. Notice that $\sigma_{\beta,3/4}(s)\equiv 1$ and $\deg H(\beta,-3/4)=1/3$.

Moreover, in Theorem \ref{thm:maass} we obtain a related result for 
the integral of $G_{\beta,m}(Z)$ against any bounded eigenfunction 
of the Laplacian on $X_L$. We intend to use this property for certain 
height pairings on the Heegner class group (see \cite{Kue3}).

In section \ref{sect4} we consider the special cases of the Siegel
 modular threefold and Hilbert modular surfaces as examples. We recover the explicit examples of \cite{Ku2} and \cite{BrKue}.

Let us finally indicate how our results can be used for 
the study of arithmetic degrees and Faltings heights. For simplicity 
(ignoring all serious technical difficulties), we assume that there exists a regular model $\mathcal{X}$ of a 
smooth compactification of $X_L$ and
a line bundle $\mathcal{M}_k$ on $\mathcal{X}$ extending the line bundle of modular
forms of weight $k$. 
Then the Hermitian line bundle $\overline{\mathcal{M}}_k$ 
given by $\mathcal{M}_k$ equipped with the (logarithmically singular) Petersson metric defines a class in a suitable arithmetic Chow ring, see \cite{BKK}.  
Using the extension of arithmetic intersection theory developed in \cite{BKK} and techniques of \cite{BrKue}, it can be shown that the 
arithmetic degree  of $\mathcal{X}$  and Faltings
  heights of subvarieties with respect to
 $\overline{\mathcal{M}}_k$
are, up to a certain rational linear combination of
logarithms of primes,  given
by linear combinations of integrals from our theorem.
The  support of these primes can be controlled by means of functoriality of the arithmetic intersection numbers 
and density results on the existence of ``many'' Borcherds products.
Details will be given in a subsequent paper.

For example, if $X_L$ is the Siegel modular variety of genus $2$, then our work
strongly  supports the following formula for the arithmetic degree:
\begin{align*}
\overline{\mathcal{M}}_k^4 
\stackrel{?}{=} k^4 \zeta(-3) \zeta(-1) 
\left( 
2\frac{\zeta'(-3)}{\zeta(-3)}+ 2\frac{\zeta'(-1)}{\zeta(-1)} + \frac{17}{6}
\right). 
\end{align*}
It can be deduced that
 the Faltings heights of Humbert surfaces $H(D) \subset
X_L$ of prime discriminant should be given by
\begin{align*}
h_{\overline{\mathcal{M}}_1}(H(D))  \stackrel{?}{=}  
 \frac{\zeta_K(-1)}{2}
\left(\frac{ \zeta_K'(-1)}{\zeta_K(-1)} +
 \frac{\zeta'(-1)}{\zeta(-1)} +\frac{3}{2}  + \frac{1}{2} \log(D) \right).
\end{align*}
In particular, our results provide further evidence for the formulas expected in \cite{MR}, \cite{Koe}, and \cite{Ku2}.
Taking into account that Humbert surfaces are birational to symmetric Hilbert modular surfaces, the latter formula essentially follows from the main theorem in \cite{BrKue}.

The idea to consider Borcherds products to
generalize results of \cite{Kue2} to higher dimension was communicated
 to the second author by S.~Kudla at the arithmetic geometry conference at 
the Issac Newton Institute in 1998. 
In Summer 1999 S.~Kudla lectured at Orsay on  
the method and  
results of \cite{Ku2} and its relation to Eisenstein series.
We would like to thank him for constantly sharing his ideas.

The approach to use the Green's functions $\Phi_{\beta,m}(Z,s)$ and results of \cite{Br1,Br2} together with the Kronecker limit formula argument arose at the annual DMV-meeting in Mainz in 1999. It can be viewed as an extension of the methods developed in \cite{Kue2}. The main work in the special case of Hilbert modular surfaces was completed at the beginning of 2001 while the first author was visiting the University of Wisconsin at Madison,
the calculations in the general case were done in summer 2001 and completed while the second author was visiting the University of Maryland.
We would like to thank these institutions for providing a stimulating atmosphere.
Moreover, we thank E.~Freitag, K.~K\"ohler, V.~Maillot, D.~Roessler, and M.~Tsuzuki for their help and useful comments.

\section{The Weil representation}
\label{sect1}

We denote by $\H=\{\tau\in \C;\; \Im(\tau)>0 \}$ the complex upper half plane. Throughout we will use $\tau$ as a standard variable on $\H$ and write $x$ for its real part and $y$ for its imaginary part, respectively.
If $z\in\C$, we abbreviate $e(z)=e^{2\pi i z}$. We denote by $\sqrt{z}=z^{1/2}$ the principal branch of the square root, so that $\arg(\sqrt{z})\in (-\pi/2,\pi/2]$. 
Moreover, for a non-zero real number $x$ we  write
$\sgn(x)=x/|x|$. 

If  $D$ is a discriminant (i.e.~a non-zero integer $\equiv 0,1\pmod{4}$), then we write $\chi_D$ for the Dirichlet character modulo $|D|$ given by the Kronecker symbol, $\chi_D(a)=\leg{D}{a}$. The corresponding Dirichlet series is denoted by $L(\chi_D,s)$.

Let $\Mp_2(\R)$ be the metaplectic cover of $\Sl_2(\R)$ realized as the group of pairs $(M,\phi)$, where $M=\kabcd\in\Sl_2(\R)$ and $\phi:\H\to \C$ is a holomorphic square root of the function $\tau\mapsto c\tau+d$ for $\tau\in \H$. 
The assignment 
\[
\abcd \mapsto \widetilde{\abcd} = \left(  \abcd , \sqrt{c\tau+d}\right)
\]
defines a locally isomorphic embedding of $\Sl_2(\R)$ into $\Mp_2(\R)$.

The inverse image $\Mp_2(\Z)$ of $\Sl_2(\Z)$ under the covering map is generated by 
\begin{align*}
T&= \left( \zxz{1}{1}{0}{1}, 1\right),\\
S&= \left( \zxz{0}{-1}{1}{0}, \sqrt{\tau}\right).
\end{align*}
There are the relations $S^2=(ST)^3=Z$, where 
$Z=\left( \kzxz{-1}{0}{0}{-1}, i\right)$
is the standard generator of the center of $\Mp_2(\Z)$.
Throughout we will frequently use the abbreviations $\Gamma_1=\Sl_2(\Z)$, $\Gamma_\infty=\{\kzxz{1}{n}{0}{1};\; n\in \Z\}\leq \Gamma_1$, and $\tilde \Gamma_\infty= \{(\kzxz{1}{n}{0}{1},1);\; n\in \Z\}\leq \Mp_2(\Z)$.

\bigskip

Let $(V,q)$ be a non-degenerate real quadratic space of signature $(b^+,b^-)$ and rank $r=b^++b^-$. We denote by $(\cdot, \cdot)$ the bilinear form corresponding to the quadratic form $q$ such that $q(x)=\frac{1}{2}(x,x)$. Let $L\subset V$ be an even lattice and write  $L'$ for its dual. Then $q$ induces a $\Q/\Z$-valued quadratic form on the finite discriminant group $L'/L$.
 
Let $(\frake_\gamma)_{\gamma\in L'/L}$ be the standard basis of the group ring $\C[L'/L]$, and write $\langle \cdot,\cdot\rangle$ for the standard scalar product on $\C[L'/L]$, which is linear in the first variable and anti-linear in the second. 
Recall that there exists a unitary representation $\rho_L$ of $\Mp_2(\Z)$ on  $\C[L'/L]$, given by
\begin{align}\label{actionT}
\rho_L(T)\frake_\gamma &= e(q(\gamma)) \frake_\gamma,\\
\label{actionS}
\rho_L(S)\frake_\gamma &= \frac{ \sqrt{i}^{b^- -b^+}}{\sqrt{|L'/L|}} \sum_{\delta\in L'/L} e(-(\gamma,\delta)) \frake_\delta
\end{align} 
(cf.~\cite{Bo2}). This representation is essentially the Weil representation associated with the quadratic module $(L'/L,q)$. 
It factors through a finite quotient of $\Mp_2(\Z)$.
Observe that 
$\rho_L(Z)\frake_\gamma= i^{b^- -b^+} \frake_{-\gamma}$.
For $\beta,\gamma\in L'/L$ and $(M,\phi)\in\Mp_2(\Z)$ we define the coefficient $\rho_{\beta\gamma}(M,\phi)$ \label{bi3} of $\rho_L$ by
\[
\rho_{\beta\gamma}(M,\phi)= \langle \rho_L(M,\phi) \frake_\gamma, \frake_\beta \rangle.
\]

Let $\kappa\in \frac{1}{2}\Z$.
The group $\Mp_2(\R)$ acts on functions $f:\H\to\C[L'/L]$ via the Petersson slash operator $f\mapsto f\mid_\kappa^* (M,\phi)$, where
\[
\big(f\mid_\kappa^* (M,\phi)\big)(\tau) = \phi(\tau)^{-2\kappa} \rho_L^*(M,\phi)^{-1} f(M\tau)
\]
and $\rho_L^*$ denotes the dual representation of $\rho_L$. 
If $f:\H\to\C[L'/L]$ is a smooth function, which is invariant under $T$, then $f$ has a Fourier expansion of the form
\[
f(\tau)=\sum_{\gamma\in L'/L}\sum_{n\in \Z-q(\gamma)} c(\gamma,n,y)\frake_\gamma(nx).
\]
Here we have abbreviated $\frake_\gamma(\tau)=\frake_\gamma e(\tau)$.

\section{Eisenstein series}
\label{sect2}

In this section we study real analytic vector valued Eisenstein series for $\Mp_2(\Z)$ transforming with the Weil representation. We compute their Fourier expansion by modifying the argument of \cite{BK}. 

Throughout we assume that $\kappa\in \frac{1}{2}\Z$, $\kappa\geq 2$,  and $2\kappa -b^- +b^+\equiv 0 \pmod{4}$. (The case $2\kappa -b^- +b^+\equiv 2 \pmod{4}$ can be treated similarly.)
If $\beta\in L'/L$ with $q(\beta)\in\Z$ and $s\in \C$, 
then the function $\H\to\C[L'/L]$ given by $\frake_\beta y^s$ is invariant under the 
action of $T,Z^2 \in \Mp_2(\Z)$.
We define the Eisenstein series $E_\beta(\tau,s)$ of weight $\kappa$ by
\begin{equation}\label{DefEb}
E_\beta (\tau,s)= \frac{1}{2}\sum_{(M,\phi)\in \tilde{\Gamma}_\infty\bs \Mp_2(\Z)} (\frake_\beta y^s) \mid_\kappa^* (M,\phi).
\end{equation}
It converges normally on $\H$ for $\Re(s)>1-\kappa/2$ and defines a $\Mp_2(\Z)$-invariant real analytic function, which is an eigenfunction of the hyperbolic Laplacian in weight $\kappa$.
Similar Eisenstein series are considered by Kudla and Yang from the adelic point of view in \cite{KuYa}.

Let $W_{\nu, \mu}(z)$ be the usual $W$-Whittaker function as in \cite{AS} Chapter 13 p.~190. For brevity we put for $s\in\C$ and $y\in \R\setminus\{0\}$:
\begin{equation}\label{eq:W}
\calW_s(y)=|y|^{-\kappa/2} W_{\sgn(y)\kappa/2,(1-\kappa)/2-s}(|y|).
\end{equation}
Notice that 
\begin{equation}\label{eq:spec}
\calW_0(y)= \begin{cases}
e^{-y/2},& \text{if $y>0$,}\\ 
e^{-y/2}\Gamma(1-\kappa,|y|),& \text{if $y<0$,}
\end{cases}
\end{equation}
where $\Gamma(a,x)=\int_x^\infty e^{-t} t^{a-1}dt$ denotes the incomplete Gamma function as in \cite{AS} p.~81.
\begin{proposition}\label{FourierE}
The Eisenstein  series $E_\beta$ has the Fourier expansion
\begin{align*}
E_\beta (\tau,s) &= \sum_{\gamma\in L'/L}\sum_{n\in \Z-q(\gamma)} c_\beta (\gamma,n,s,y)  \frake_\gamma(n x),
\end{align*}
where the coefficients  $c_\beta (\gamma,n,s,y)$ are equal to 
\begin{equation*}
\begin{cases}
\displaystyle
(\delta_{\beta,\gamma}+\delta_{-\beta,\gamma})y^s
+2 \pi  y^{1-\kappa-s}\frac{\Gamma(\kappa+2s-1)}{\Gamma(\kappa+s)\Gamma(s)}\sum_{c\in\Z-\{0\}}|2c|^{1-\kappa-2s} H_c^*(\beta,0,\gamma,0),& n=0,\\
\displaystyle
\frac{2^\kappa \pi^{s+\kappa}  |n|^{s+\kappa-1}}{\Gamma(s+\kappa)}\calW_s(4\pi ny)\sum_{c\in\Z-\{0\}}|c|^{1-\kappa-2s} H_c^*(\beta,0,\gamma,n), &n>0,\\
\displaystyle
\frac{2^\kappa \pi^{s+\kappa}  |n|^{s+\kappa-1}}{\Gamma(s)}\calW_s(4\pi ny) \sum_{c\in\Z-\{0\}}|c|^{1-\kappa-2s} H_c^*(\beta,0,\gamma,n), &n<0,
\end{cases}
\end{equation*}
and $H_c^*(\beta,m,\gamma,n)$ denotes the generalized Kloosterman sum
\begin{equation}\label{DefH_c*}
H_c^*(\beta,m,\gamma,n)=\frac{e^{-\pi i\sgn(c)\kappa/2}}{|c|} \sum_{\substack{d\;(c)^*\\ \kabcd\in\Gamma_\infty\bs \Gamma_1/\Gamma_\infty}} \rho_{\beta\gamma}\widetilde{\abcd}\, e\!\left( \frac{ma+nd}{c}\right).
\end{equation}
\end{proposition}

The sum in (\ref{DefH_c*}) runs over all primitive residues $d$ modulo $c$ and $\kabcd$ is a representative for the double coset in $\Gamma_\infty\bs \Gamma_1/\Gamma_\infty$ with lower row $(c\; d')$ and $d'\equiv d\pmod{c}$. Observe that the quantity $\rho_{\beta\gamma}\widetilde{\kabcd}\, e\!\left( \tfrac{ma+nd}{c}\right)$ does not depend on the choice of the coset representative.

The coefficients $\rho_{\beta\gamma}\widetilde{\kabcd}$ are universally bounded, since $\rho_L$ factors through a finite group. Hence there is a constant $C>0$ such that $H_c^*(\beta,0,\gamma,n)<C$ for all $\gamma \in L'/L$, $n\in\Z-q(\gamma)$, and $c\in\Z-\{0\}$. This implies that the series for the Fourier coefficients converge absolutely for $\Re(s)>1-\kappa/2$.

\begin{proof}
Let $\gamma\in L'/L$ and $n\in\Z-q(\gamma)$.
We argue as in \cite{Br2} Chapter 1.2.3. The coefficient $c_\beta(\gamma,n,s,y)$ is given by the Fourier integral
\[
c_\beta(\gamma,n,s,y)=\frac{1}{2}\int\limits_0^1  \sum_{(M,\phi)\in \tilde{\Gamma}_\infty\bs \Mp_2(\Z)} \big\langle (\frake_\beta y^s) \mid_\kappa^* (M,\phi),\; \frake_\gamma(nx)\big\rangle\,dx.
\]
We split the above sum into the sum over $1,Z,Z^2,Z^3\in\tilde{\Gamma}_\infty\bs \Mp_2(\Z)$ and the sum over $(\kabcd,\phi)\in\tilde{\Gamma}_\infty\bs \Mp_2(\Z)$ with $c\neq 0$. 
Since $\frake_\beta y^s$ is invariant under the action of $Z^2$, we get
\begin{align*}
 c_\beta(\gamma,n,s,y)
&=\delta_{0,n}(\delta_{\beta,\gamma}+\delta_{-\beta,\gamma})y^s
+ \sum_{\substack{c\neq 0\\ \kabcd\in \Gamma_\infty\bs \Gamma_1/\Gamma_\infty }}\rho_{\beta\gamma}\widetilde{\abcd}
\int\limits_{-\infty}^\infty \frac{y^s e(-nx)}{(c\tau+d)^{\kappa+s}(c\bar\tau+d)^{s}}\,dx.
\end{align*}
We now compute the integral. Because $\sqrt{c\tau+d}=\sgn(c)\sqrt{c}\sqrt{\tau+d/c}$, we have
\begin{align*}
\int\limits_{-\infty}^\infty \frac{y^s e(-nx)}{(c\tau+d)^{\kappa+s}(c\bar\tau+d)^{s}}\,dx
&= |c|^{-\kappa-2s}\sgn(c)^\kappa y^s \int\limits_{-\infty}^\infty \frac{ e(-nx)}{(\tau+d/c)^{\kappa+s}(\bar\tau+d/c)^{s}}\,dx \\
&= |c|^{-\kappa-2s}\sgn(c)^\kappa e(nd/c) y^s \int\limits_{-\infty}^\infty \frac{ e(-nx)}{\tau^{\kappa+s}\bar\tau^{s}}\,dx. 
\end{align*}
Consequently
\begin{equation}
\label{eq1}
 c_\beta(\gamma,n,s,y)
=\delta_{0,n}(\delta_{\beta,\gamma}+\delta_{-\beta,\gamma})y^s
+i^\kappa y^s \int\limits_{-\infty}^\infty \frac{ e(-nx)}{\tau^{\kappa+s}\bar\tau^{s}}\,dx
\sum_{c\in\Z-\{0\}}|c|^{1-\kappa-2s} H_c^*(\beta,0,\gamma,n).
\end{equation}
Using \cite{B2} 3.2 (12) on p. 119 (and correcting the sign of the first formula there), we find for the latter integral
\begin{align*}
\int\limits_{-\infty}^\infty \frac{ e(-nx)}{\tau^{\kappa+s}\bar\tau^{s}}\,dx
&=i^{-\kappa}\int\limits_{-\infty}^\infty \frac{ e^{-2\pi inx}}{(y-ix)^{\kappa+s}(y+ix)^{s}}\,dx\\
&=2\pi i^{-\kappa}\begin{cases} 
(2y)^{-\kappa/2-s}\Gamma(\kappa+s)^{-1}(2\pi n)^{s-1+\kappa/2} W_{\kappa/2, (1-\kappa)/2-s}(4\pi n y),& n>0,\\
(2y)^{-\kappa/2-s}\Gamma(s)^{-1}(2\pi |n|)^{s-1+\kappa/2} W_{-\kappa/2, (1-\kappa)/2-s}(4\pi |n| y),& n<0,
\end{cases}\\
&=2^\kappa \pi^{s+\kappa}i^{-\kappa}  |n|^{s+\kappa-1} y^{-s} \calW_s(4\pi ny) 
\begin{cases}\Gamma(\kappa+s)^{-1},& n>0,\\
\Gamma(s)^{-1},& n<0.
\end{cases}
\end{align*}
For $n=0$ we derive by means of  \cite{Fr} Chap.~III Lemma 4.4 that
\begin{align*}
\int\limits_{-\infty}^\infty \frac{ 1}{\tau^{\kappa+s}\bar\tau^{s}}\,dx
&=i^{-\kappa} y^{1-2s-\kappa}\int\limits_{-\infty}^\infty \frac{1}{(1-ix)^\kappa |1-ix|^{2s}}\,dx\\
&=2^{2-\kappa-2s} \pi i^{-\kappa} \frac{\Gamma(\kappa+2s-1)}{\Gamma(\kappa+s)\Gamma(s)} y^{1-\kappa-2s}.
\end{align*}
Inserting into \eqref{eq1}, we obtain the assertion.
\end{proof}

The Weil representation $\rho_L$ is trivial on the principal congruence subgroup $\Gamma(N)$ of level $N$ of $\Mp_2(\Z)$, where $N$ is the least positive integer such that $Nq(\gamma)\in \Z$ for all $\gamma\in L'/L$. So the Eisenstein series $E_\beta(\tau,s)$ is a finite linear combination of $\Mp_2(\Z)/\Gamma(N)$-translates of classical scalar valued Eisenstein series for $\Gamma(N)$. This implies that $E_\beta(\tau,s)$ has a meromorphic continuation in $s$ to the whole complex plane, which is holomorphic in $s$ at $s=0$. 

Because of Proposition \ref{FourierE} and \eqref{eq:spec} the special value $E_\beta(\tau,0)$ is always a holomorphic modular form of weight $\kappa$, if $\kappa>2$. For $\kappa=2$ it is holomorphic, if $\rho_L$ does not contain the trivial representation as an irreducible constituent.

\bigskip

We are mainly interested in the Eisenstein series $E_0(\tau,s)$.
Using the ideas of \cite{BK}, its coefficients can be computed more explicitly. We briefly indicate the main steps.
By means of Shintani's formula \cite{Sh} for the coefficients of the Weil representation it can be shown that
\begin{equation*}
H_c^*(0,0,\gamma,n)=\frac{(-1)^{(2\kappa-b^-+b^+)/4}}{\sqrt{|L'/L|}} |c|^{-1+r/2}\sum_{a\mid c} a^{1-r} \mu(|c|/a) N_{\gamma,n}(a),
\end{equation*}
where $\mu$ denotes the Moebius function and $N_{\gamma,n}(a)$ the representation number 
\begin{equation}\label{na}
N_{\gamma,n}(a)=\# \{ x\in L/aL;\quad q(x-\gamma)+n\equiv 0\pmod{a}\}
\end{equation}
modulo $a$ (see \cite{BK} Proposition 3).
As a consequence we may infer that
\begin{equation*}
\sum_{c\in\Z-\{0\}}|c|^{1-\kappa-2s} H_c^*(0,0,\gamma,n)
= \frac{2(-1)^{(2\kappa-b^-+b^+)/4}}{\sqrt{|L'/L|}\zeta(2s+\kappa-r/2)}
\sum_{a=1}^\infty N_{\gamma,n}(a)a^{1-r/2-\kappa-2s}
\end{equation*}
(see \cite{BK} Proposition 4). Here $\zeta(s)$ is the Riemann zeta function.

We now compute the $L$-series
\begin{equation}
\label{step2}
L_{\gamma,n}(s)=\sum_{a=1}^\infty N_{\gamma,n}(a)a^{-s}.
\end{equation}
Since the representation number function $N_{\gamma,n}(a)$ is multiplicative 
in $a$, this $L$-series has an Euler product expansion. 
We let 
\begin{equation*}
d_\gamma=\min\{b\in \Z_{>0};\quad b\gamma\in L\}
\end{equation*} 
be the level of $\gamma$. Then $d_\gamma$ clearly divides $|L'/L|$, and $d_\gamma^2 n$ and $2d_\gamma n$ are integers.
If we put  
\begin{equation}\label{wp}
w_p=1+2v_p(2nd_\gamma),
\end{equation}
where $v_p$ denotes the (additive) 
$p$-adic valuation on $\Q$,
then $N_{\gamma,n}(p^{\alpha+1}) = p^{r-1}N_{\gamma,n}(p^\alpha)$
for all $\alpha\geq w_p$.
We may conclude that
\[
L_{\gamma,n}(s)=\zeta(s-r+1)\prod_{\text{$p$ prime}} L_{\gamma,n}^{(p)}(p^{-s}),
\]
where $L_{\gamma,n}^{(p)}(X)$ denotes the polynomial
\begin{equation}\label{eq:locl}
L_{\gamma,n}^{(p)}(X)=N_{\gamma,n}(p^{w_p})X^{w_p}+(1-p^{r-1}X)\sum_{\nu=0}^{w_p-1}N_{\gamma,n}(p^\nu)X^\nu
\in \Z[X].
\end{equation}
Observe that $L_{\gamma,n}^{(p)}(X)$ does not change, if $w_p$ is replaced by any greater integer.

\begin{proposition}\label{prop:eisex}
The coefficient $c_0(\gamma,n,s,y)$ of $E_0(\tau,s)$ is equal to 
\[
\begin{cases}
\displaystyle
2 \delta_{0,\gamma}y^s 
{}+2^{3-\kappa-2s} \pi  y^{1-\kappa-s}\frac{\Gamma(\kappa+2s-1)}{\Gamma(\kappa+s)\Gamma(s)}
\frac{(-1)^{(2\kappa-b^-+b^+)/4}}{\sqrt{|L'/L|}}
\prod_{p} L^{(p)}_{\gamma,n}(p^{1-r/2-\kappa-2s}),& n=0,\\
\displaystyle
\frac{2^{\kappa+1} \pi^{s+\kappa}  |n|^{s+\kappa-1}}{\Gamma(s+\kappa)} \calW_s(4\pi ny)
\frac{(-1)^{(2\kappa-b^-+b^+)/4}}{\sqrt{|L'/L|}}
\prod_{p} L^{(p)}_{\gamma,n}(p^{1-r/2-\kappa-2s})
, &n>0,\\
\displaystyle
\frac{2^{\kappa+1} \pi^{s+\kappa}  |n|^{s+\kappa-1}}{\Gamma(s)} \calW_s(4\pi ny)
\frac{(-1)^{(2\kappa-b^-+b^+)/4}}{\sqrt{|L'/L|}}
\prod_{p} L^{(p)}_{\gamma,n}(p^{1-r/2-\kappa-2s})
, &n<0.
\end{cases}
\]
\end{proposition} 

For $n\neq 0$ we define coefficients $C(\gamma,n,s)$ by
\begin{equation}\label{def:C}
c_0(\gamma,n,s,y)=C(\gamma,n,s)\calW_s(4\pi ny),
\end{equation}
so that $C(\gamma,n,s)$ is essentially the part of $c_0(\gamma,n,s,y)$ coming from the finite places of $\Q$.
In this case the Euler products in Proposition \ref{prop:eisex} 
can be computed 
using a result of Siegel on representation numbers of quadratic forms modulo prime powers (\cite{Si} Hilfssatz 16). One finds that if $p$ is a prime coprime to the integer $2d_\gamma^2 n |L'/L|$ (in particular $w_p=1$), then
\[
L_{\gamma,n}^{(p)}(X)=
\begin{cases} \displaystyle 1 -\chi_D(p)p^{r/2-1}X ,& \text{if $r$ is even,}\\[1ex]
 \displaystyle \frac{1-p^{r-1}X^2}{1-\chi_D(p)p^{r/2-1/2}X} ,& \text{if $r$ is odd.}\end{cases}
\]
Here $D$ is the discriminant
\[
\begin{cases} 
\displaystyle (-1)^{r/2}\det(L) ,& \text{if $r$ is even,}\\[1ex]
 \displaystyle 2(-1)^{(r+1)/2}d_\gamma^2 n \det(L),& \text{if $ r$ is odd,}
\end{cases}
\]
and $\det(L)$ means the Gram determinant of $L$. 
We may write $D$ uniquely as the product $D=D_0f^2$ of a fundamental discriminant $D_0$ and the square of a positive integer $f$. Then $D_0$ is the discriminant of the quadratic number field $\Q(\sqrt{D})$, and the character $\chi_{D_0}$ is primitive.
We find that, up to a finite product of rational functions in $p^{-s}$ over primes $p$ dividing $2d_\gamma^2 n |L'/L|$, the Euler product $\prod_p L^{(p)}_{\gamma,n}(p^{1-r/2-\kappa-2s})$ is equal to
\begin{equation}\label{eq:euler}
\begin{cases} 
\displaystyle \frac{1}{L(\chi_{D_0},2s+\kappa)} ,& \text{if $r$ is even,}\\[3ex]
 \displaystyle \frac{L(\chi_{D_0},2s+\kappa-1/2)}{\zeta(4s+2\kappa-1)} ,& \text{if $ r$ is odd.}
\end{cases}
\end{equation}

\begin{theorem}\label{thm:eismain}
If $n>0$, then $
C(\gamma,n,s)$ is equal to
\begin{align*}
&2^{2-2s} \pi^{-s}  |n|^{\kappa+s-1} 
\cos\big(\pi(s-\tfrac{\delta}{2}+\tfrac{b^--b^+}{4})\big)
\frac{|D_0|^{2s+\kappa-1/2} \Gamma(2s+\kappa)}{\sqrt{|L'/L|}\Gamma(s+\kappa)}
\cdot \frac{\sigma_{\gamma,n}(2s+\kappa)}{L(\chi_{D_0},1-2s-\kappa)},  
\end{align*}
if $r$ is even, and to
\begin{align*}
&2^{2s+2\kappa-1/2}\pi^{-s}|n|^{\kappa+s-1} 
\\
&\times
\frac{\sin\big(\pi(2s+\kappa)\big)|D_0|^{1-\kappa-2s}\Gamma(2s+\kappa)}{\cos\big(\pi(s-\tfrac{\delta}{2}+\tfrac{b^--b^+-1}{4})\big)\sqrt{|L'/L|}\Gamma(s+\kappa)}
\cdot\frac{L(\chi_{D_0},3/2-2s-\kappa)}{\zeta(2-4s-2\kappa)} \sigma_{\gamma,n}(2s+\kappa),   
\end{align*}
if $r$ is odd. 
Here $D_0$ denotes the fundamental discriminant defined above and $\delta=0$, if $D_0>0$, and $\delta=1$, if $D_0<0$.
Moreover, $\sigma_{\gamma,n}(s)$ is the generalized divisor sum 
\begin{equation}\label{eq:defsigma}
\sigma_{\gamma,n}(s)=
\begin{cases} 
\displaystyle \prod_{p\mid 2d_\gamma^2 n \det(L)}\frac{L^{(p)}_{\gamma,n}(p^{1-r/2-s})}{1-\chi_{D_0}(p)p^{-s}} ,& \text{if $r$ is even,}\\[3ex]
 \displaystyle  \prod_{p\mid 2d_\gamma^2 n \det(L)}
\frac{1-\chi_{D_0}(p)p^{1/2-s}}{1-p^{1-2s}}L^{(p)}_{\gamma,n}(p^{1-r/2-s}) ,& \text{if $ r$ is odd.}
\end{cases}
\end{equation}
Analogous formulas hold for $n<0$, but with $\Gamma(s+\kappa)$ replaced by $\Gamma(s)$.
\end{theorem}

\begin{proof}

From Proposition \ref{prop:eisex} and the above discussion we may infer that $C(\gamma,n,s)$ is equal to
\begin{equation}\label{eq:for}
\frac{(-1)^{(2\kappa-b^-+b^+)/4}2^{\kappa+1} \pi^{s+\kappa}  |n|^{s+\kappa-1}}{\sqrt{|L'/L|}\Gamma(s+\kappa)}
\cdot
\begin{cases}
\displaystyle \frac{1}{L(\chi_{D_0},2s+\kappa)}\sigma_{\gamma,n}(2s+\kappa) ,& \text{if $r$ is even,}\\[3ex]
 \displaystyle \frac{L(\chi_{D_0},2s+\kappa-1/2)}{\zeta(4s+2\kappa-1)}\sigma_{\gamma,n}(2s+\kappa),& \text{if $ r$ is odd.}
\end{cases}
\end{equation}
Because the character $\chi_{D_0}$ is primitive, the Dirichlet series $L(\chi_{D_0},s)$ satisfies the functional equation
\[
L(\chi_{D_0},s)=\frac{2^{s-1} \pi^s |D_0|^{1/2-s}}{\cos(\pi (s-\delta)/2)\Gamma(s)}L(\chi_{D_0},1-s).
\]
In fact, this form of the functional equation can be deduced from the standard form (see e.g.~\cite{Za} p.~53) by means of the identity $\Gamma(s)\Gamma(1-s)=\pi/\sin(\pi s)$ and the Legendre duplication formula.
In particular, the Riemann zeta function has the functional equation
\[
\zeta(s)=\frac{2^{s-1} \pi^s}{\cos(\pi s/2)\Gamma(s)}\zeta(1-s).
\]
If we insert this into \eqref{eq:for}, we see that $C(\gamma,n,s)$ is equal to
\begin{align*}
&(-1)^{(2\kappa-b^-+b^+)/4}2^{2-2s} \pi^{-s}  |n|^{\kappa+s-1} \frac{|D_0|^{2s+\kappa-1/2}}{\sqrt{|L'/L|}}\\
&\times\cos\big(\pi(s+\kappa/2-\delta/2)\big)
\frac{\Gamma(2s+\kappa)}{\Gamma(s+\kappa)}
\frac{\sigma_{\gamma,n}(2s+\kappa)}{L(\chi_{D_0},1-2s-\kappa)},  
\end{align*}
if $r$ is even, and to
\begin{align*}
&(-1)^{(2\kappa-b^-+b^+)/4} 2^{2s+2\kappa-1/2}\pi^{-s}|n|^{\kappa+s-1} 
\frac{|D_0|^{1-\kappa-2s}}{\sqrt{|L'/L|}}\\
&\times
\frac{\cos\big(\pi(2s+\kappa-1/2)\big)}{\cos\big(\pi(s+\kappa/2-1/4-\delta/2)\big)}
\frac{\Gamma(2s+\kappa)}{\Gamma(s+\kappa)} 
\frac{L(\chi_{D_0},3/2-2s-\kappa)}{\zeta(2-4s-2\kappa)} \sigma_{\gamma,n}(2s+\kappa),   
\end{align*}
if $r$ is odd. This easily implies the assertion.
\end{proof}

The finite Euler products $\sigma_{\gamma,n}(s)$ are interesting invariants of the lattice $L$ (or better of $L'/L$ as we shall see below). Notice that in 
their definition  one could replace the product over primes dividing  $2d_\gamma^2 n \det(L)$ by the product over any larger set of primes including the primes dividing $2d_\gamma^2 n \det(L)$. 
The part of $\sigma_{\gamma,n}(s)$ given by primes dividing $d_\gamma^2 n$ but coprime to $2\det(L)$ can also be evaluated just in terms of $n,\gamma, \det(L)$, not involving representation numbers (see \cite{BK} Theorem 7). 

We have the following table for $\sgn(D_0)$ (and hence for $\delta$):
\begin{table}[h]\label{table}
\begin{tabular}{l||c|c|c|c}
$b^+-b^- \pmod{4}$   & \; 0 \; & \; 1\; & \; 2\; & \; 3 \;\\
\hline\hline
$n>0$                & $+$ & $-$ & $-$ & $+$ \\ \hline
$n<0$                & $+$ & $+$ & $-$ & $-$   
\end{tabular}
\end{table}

The following proposition is very useful for explicit computations.

\begin{proposition} \label{prop:genus}
Suppose that  $L$ and $\tilde L$ are even lattices of rank $r$ and $\tilde r$, respectively. Write $L_{\gamma,n}^{(p)}(X)$ and $\tilde L_{\gamma,n}^{(p)}(X)$ for the $p$-polynomials defined in \eqref{eq:locl} corresponding to $L$ and $\tilde L$, respectively. If $L'/L\cong \tilde L'/\tilde L$, then
\[
L_{\gamma,n}^{(p)}(p^{-r/2} X)=\tilde L_{\gamma,n}^{(p)}(p^{-\tilde r/2} X).
\]
Consequently, if $\sigma_{\gamma,n}(s)$ and $\tilde \sigma_{\gamma,n}(s)$ denote the finite Euler products defined in \eqref{eq:defsigma} corresponding to $L$ and $\tilde L$, then 
$\sigma_{\gamma,n}(s)=\tilde \sigma_{\gamma,n}(s)$. 
\end{proposition}

\begin{proof}
This can  be deduced from the fact that the definition \eqref{DefEb} of $E_0(\tau,s)$ only depends on $L'/L$ (via the Weil representation). The quantities apart from  $\sigma_{\gamma,n}(s)$ in the formula for $C(\gamma,n,s)$ in Theorem \ref{thm:eismain} also only depend on $L'/L$. Hence the same must be true for $\sigma_{\gamma,n}(s)$. We leave it to the reader to fill in  the details.
\end{proof}

\section{Green's functions associated with Heegner divisors}
\label{sect3}

For the rest of this paper let $(V,q)$ be a real quadratic space of signature $(2,p)$ and put $\kappa=1+p/2$. 
We assume that either $p> 2$, or that $p=2$ and the dimension of a maximal isotropic subspace of $L\otimes_\Z\Q$ equals $1$. (We will use many results of \cite{Br1}, where the general assumption was $p>2$ for simplicity. However, it is easily verified that everything we need here also works for the latter case. See also \cite{Br2}.) 
In this section, beside other things, we will relate the values of integrals over certain Green's functions to the derivatives of the coefficients of the Eisenstein series $E_0(\tau,s)$ in weight $\kappa$.
  
We start by recalling some facts on the Hermitean symmetric space associated with the real orthogonal group $\Orth(V)$ of  $(V,q)$.
It can be realized as the Grassmannian
$\Gr(V)$ of oriented positive definite subspaces $v\subset V$ of dimension $2$.
It is a real analytic manifold of real dimension $2p$ consisting of $2$ connected components given by the $2$ possible choices of the orientation.
The complex structure on $\Gr(V)$ can be most easily realized as follows. We consider the complexification $V_\C= V\otimes_\R\C$ of $V$ and extend the bilinear form $(\cdot,\cdot)$  to a $\C$-bilinear form on $V_\C$. The subset
\[
\calK=\{ [Z]\in P(V_\C); \quad (Z,Z)=0, \; (Z,\bar Z)>0\}
\]
of the projective space over $V_\C$ is invariant under the action of $\Orth(V)$. By mapping a positively oriented orthogonal basis $X,Y$ of $v\in \Gr(V)$ with $X^2=Y^2>0$ to $[X+iY]\in P(V_\C)$, one obtains a real analytic isomorphism $\Gr(V)\to \calK$. It defines a complex structure on the Grassmannian.

We fix a connected component $\Gr'(V)$ of $\Gr(V)$ and write $\calK'$ for the corresponding component of $\calK$.
A subgroup of $\Orth(V)$ of index $2$, denoted by $\Orth'(V)$, preserves these components. It is readily verified that $\Orth'(V)$ is the set of elements whose spinor norm has the same sign as the determinant.
Let $\Orth'(L)=\Orth(L)\cap\Orth'(V)$, and write $\Gamma(L)$ for the discriminant kernel of $\Orth'(L)$, i.e.,  the kernel of the natural homomorphism $\Orth'(L)\to\Orth(L'/L)$.
By the theory of Baily-Borel, the quotient 
\[
X_L=\Gamma(L)\bs \calK'
\]
is a quasi-projective algebraic variety over $\C$ of dimension $p$. 
     
There are certain special divisors on $X_L$ arising from embedded quotients analogous to $X_L$ of dimension $p-1$. 
For any  vector $\lambda\in V$ of negative norm the orthogonal complement $\lambda^\perp$ in $\Gr(V)$ defines a complex analytic divisor on $\Gr(V)$.
If $\beta\in L'/L$, and $m\in \Z+q(\beta)$ is negative, then 
\[
H(\beta,m)=\sum_{\substack{\lambda\in L+\beta\\q(\lambda)=m}}\lambda^\perp
\]
is a $\Gamma(L)$-invariant divisor on $\Gr(V)$ called the {\em Heegner divisor} of discriminant $(\beta,m)$. Its restriction to $\Gr'(V)$ is the inverse image under the canonical projection of an algebraic divisor on $X_L$, also denoted by $H(\beta,m)$. Notice that $H(\beta,m)=H(-\beta,m)$, and that the multiplicities of all irreducible components of $H(\beta,m)$ are $2$, if $2\beta=0$, and $1$, if $2\beta\neq 0$ in $L'/L$.

We consider automorphic Green's functions for these Heegner divisors. They 
were first introduced in \cite{Br1,Br2}, and independently from a different perspective in \cite{OT}.  
We define the Green's function associated with $H(\beta,m)$  for $v\in \Gr(V)\setminus H(\beta,m)$ and $s\in\C$ with $\Re(s)>\kappa/2$ by 
\begin{align}
\nonumber
\Phi_{\beta,m}(v,s)  = 2\frac{\Gamma(s-1+\kappa/2)}{\Gamma(2s)} &\sum_{\substack{\lambda\in \beta+L\\ q(\lambda)=m}}  \left(\frac{m}{m-q(\lambda_v)}\right)^{s-1+\kappa/2}\\
\label{def:green}
&\phantom{=}
\times F\left(s-1+\kappa/2,\, s+1-\kappa/2, \, 2s; \, \frac{m}{m-q(\lambda_v)}\right).
\end{align}
Here $\lambda_v$ denotes the orthogonal projection of $\lambda$ to $v\subset V$ and 
\begin{align}\label{def:gauss}
F(a,b,c; z) &= \sum_{n=0}^\infty \frac{(a)_n (b)_n}{(c)_n} \frac{z^n}{n!}
\end{align}
(with $(a)_n = \Gamma(a+n)/\Gamma(a)$) the Gauss hypergeometric function as in 
\cite{AS} Chap.~15.
In \cite{Br1} the function $\Phi_{\beta,m}(v,s)$ was originally defined as a regularized theta lifting of a certain non-holomorphic Poincar\'e series of weight $k=2-\kappa$ and then shown to be equal to the series above. (See \cite{Br1} Theorem 2.14. Note that the regularization due to Borcherds and Harvey-Moore is replaced by a certain {\em automorphic} regularization.)
Here we prefer to take \eqref{def:green} as our definition, since it is a bit more direct. This is also the approach of 
\cite{Br2} and \cite{OT}.
The series converges normally for $v\in \Gr(L)\setminus H(\beta,m)$ and $\Re(s)>\kappa/2$ and defines a $\Gamma(L)$-invariant function. It can be easily written in the coordinates of $\calK$ by noting that $q(\lambda_v)=|(\lambda,Z)|^2/(Z,\bar Z)$ for $[Z]\in \calK$ corresponding to $v$. 
We denote the resulting  $\Gamma(L)$-invariant function on  $\calK$ by $\Phi_{\beta,m}(Z,s)$. 

According to \cite{Br1} Prop.~2.8 the function $\Phi_{\beta,m}(v,s)$ has a 
meromorphic continuation in $s$ to a neighborhood of $\kappa/2$. It has a simple pole at $s=\kappa/2$.
As a function in $v$, it is real analytic on $\Gr(V)\setminus H(\beta,m)$ and has a logarithmic singularity along $H(\beta,m)$. If $U\subset\Gr(V)$ denotes a compact neighborhood of $v$, then there is a finite set $S(U)$ of $\lambda\in \beta +L$ with $q(\lambda)=m$ such that 
\[
\Phi_{\beta,m}(v,s)= -4 \sum_{\lambda\in S(U)} \log|\lambda_v| +O(1)
\]
on $U$ (cf.~\cite{Br1} Theorem 2.12 and \cite{Bo2} Theorem 6.2).

Let $\Delta$ be the $\Orth(V)$-invariant Laplace operator on $\Gr(V)$ induced by the Casimir element of the Lie algebra of $\Orth(V)$. We normalize it as in \cite{Br1}. 
The Green's function $\Phi_{\beta,m}(v,s)$ is an eigenfunction of $\Delta$. More precisely, according to \cite{Br1} Theorem 4.6, we have
\begin{align}\label{eq:lapl}
\Delta \Phi_{\beta,m}(v,s) = \frac{1}{2}(s-\kappa/2)(s-1+\kappa/2)\Phi_{\beta,m}(v,s).
\end{align}

The Fourier expansion of $\Phi_{\beta,m}(v,s)$ can be computed using the integral representation as a regularized theta integral. 
Let $z\in L$ be a primitive isotropic vector and $N$ be the unique positive integer such that $(z,L)=N\Z$. Let $z'\in L'$ with $(z,z')=1$.
Then $K=L\cap z^\perp\cap {z'}^\perp$ is an even lattice of signature $(1,p-1)$ isomorphic to $(L\cap z^\perp)/\Z z$.
The subset 
\[
\calH_z=\{ Z\in K\otimes_\Z \C; \quad q(\Im( Z))>0\}
\]
of $K\otimes_\Z \C\cong\C^p$, is a tube domain realization of the Hermitean symmetric space $\calK$, the isomorphism being given by 
$Z\mapsto [Z+z'-q(Z)z-q(z')z]$. The invariance of a function $F$ on $\calK$ under certain Eichler transformations in $\Gamma(L)$ implies that the corresponding function $F_z$ on $\calH_z$ is periodic with period lattice $K'$, i.e.,
$ F_z(Z+\lambda)=F_z(Z)$ for all $\lambda\in K'$. In particular any sufficiently smooth $\Gamma(L)$-invariant function on $\calK$ has a Fourier expansion with Fourier coefficients indexed by the lattice $K'$. We call it the Fourier expansion at the cusp $z$. The coefficient corresponding to $\lambda=0\in K'$ is loosely called the constant Fourier coefficient at the cusp $z$. 
 
According to \cite{Br1} Theorem 2.15 and Lemma 1.15 the constant coefficient at the cusp $z$ of the function $\Phi_{\beta,m}(Z,s)$ is given by
\begin{align}\label{const}
 &  U_{\beta,m}(Y,s)+2(Y^2)^{\kappa/2-s}\varphi_{\beta,m}(s)\zeta(2s+1-\kappa),
\end{align}
where
\begin{align}
\label{eq:varphi}
\varphi_{\beta,m}(s)&=\frac{1}{\sqrt{\pi}}\left( \frac{2}{\pi }\right)^{s-\kappa/2} \Gamma(s+1/2-\kappa/2)  b( 0, 0,s),\\
\label{def:U}
U_{\beta,m}(Y,s)&=\frac{|Y|}{\sqrt{2} } \Phi^K_{\beta,m}(Y/|Y|,s)\\
\nonumber
&\phantom{=}{}+ \frac{2}{\sqrt{\pi}}\left( \frac{2}{\pi Y^2 }\right)^{s-\kappa/2} \Gamma(s+1/2-\kappa/2) \sum_{\ell=1}^{N-1} b( \ell z/N, 0,s) \sum_{n\geq 1}e(\ell n/N)n^{\kappa-1-2s},\\
\label{def:b}
b(\gamma,0,s) &= -\frac{2^\kappa\pi^{s+\kappa/2} |m|^{s-1+\kappa/2}}{(2s-1)\Gamma(s+1-\kappa/2)\Gamma(s-1+\kappa/2)} \sum_{c\in\Z-\{0\}} |c|^{1-2s} H_{-c}^*(\gamma,0,\beta,-m),
\end{align}
and $Z\in\calH_z$ with $Y=\Im(Z)$.
Here we have used that $k=2-\kappa$, $|z_v|=1/|Y|$, and $w=Y/|Y|$.
The function 
$U_{\beta,m}(Y,s)$ is continuous in $Y$ and holomorphic in $s$ at $\kappa/2$.  The second summand in \eqref{const} will be of particular interest to us. It has a simple pole with residue $\varphi_{\beta,m}(\kappa/2)$ at $s=\kappa/2$. Observe that $\varphi_{\beta,m}(s)$ is independent of the choice of $z$.
The difference of $\Phi_{\beta,m}(Z,s)$ and its constant coefficient at $z$ is holomorphic in $s$ near $\kappa/2$.

We will be interested in the following regularized Green's functions.

\begin{definition}
Let $\beta\in L'/L$ and $m\in \Z+q(\beta)$ with $m<0$. We define
the regularized Green's function $\Phi_{\beta,m}(v)$ associated with the Heegner divisor $H(\beta,m)$ to be the constant term in the Laurent expansion in $s$ of $\Phi_{\beta,m}(v,s)$ at $s=\kappa/2$.
\end{definition}

\begin{proposition}\label{prop:aux}
The regularized Green's function $\Phi_{\beta,m}(v)$ is given by
\[
\Phi_{\beta,m}(v)=\lim_{s\to \kappa/2}\left(\Phi_{\beta,m}(v,s)-\frac{\varphi_{\beta,m}(\kappa/2)}{s-\kappa/2}\right).
\]
Its constant Fourier coefficient at the cusp $z$ is equal to
\[
U_{\beta,m}(Y,\kappa/2)-\varphi_{\beta,m}(\kappa/2)\log(Y^2) +\varphi'_{\beta,m}(\kappa/2)-2\varphi_{\beta,m}(\kappa/2)\Gamma'(1).
\]
\end{proposition}

\begin{proof}
To compute $\Phi_{\beta,m}(v)$ we need the Laurent expansion of the second summand \[
2(Y^2)^{\kappa/2-s}\varphi_{\beta,m}(s)\zeta(2s+1-\kappa)
\]
of \eqref{const} at $s=\kappa/2$.
Using the Laurent expansion $\zeta(s)=1/(s-1)-\Gamma'(1)+O(s-1)$ of the Riemann zeta function at $s=1$, we see that it is given by
\[
\frac{\varphi_{\beta,m}(\kappa/2)}{s-\kappa/2}
- \varphi_{\beta,m}(\kappa/2)\log(Y^2)
+ \varphi'_{\beta,m}(\kappa/2)-2\varphi_{\beta,m}(\kappa/2)\Gamma'(1)+O(s-\kappa/2).
\] 
In view of the properties of the constant Fourier coefficient  of $\Phi_{\beta,m}(v,s)$ at the cusp $z$, this implies the assertion.
\end{proof}

To ease the comparison we notice that the constant $C_{\beta,m}$  in \cite{Br1} Theorem 3.9 (which depends on the choice of the cusp $z$) is given by
\begin{equation}\label{eq:Cb}
C_{\beta,m}=\varphi'_{\beta,m}(\kappa/2)-2\varphi_{\beta,m}(\kappa/2)\Gamma'(1)-2\sum_{\ell=1}^{N-1} b(\ell z/N,0,\kappa/2)\log |1-e(\ell/N)|.
\end{equation}

\bigskip

We now relate $\varphi_{\beta,m}(s)$ to the coefficients of the Eisenstein series $E_0(\tau,s)$ in weight $\kappa=1+p/2$ of the previous section. 
If $\Phi_{\beta,m}(v,s)$ is viewed as the regularized theta lifting of a certain non-holomorphic Poincar\'e series $F_{\beta,m}(\tau,s)$ for $\Mp_2(\Z)$ as in \cite{Br1}, then $b(0,0,s)$ naturally arises as the constant Fourier coefficient of $F_{\beta,m}(\tau,s)$.
The point is that this constant coefficient is up to a universal factor equal to $C(\beta,-m,s-\kappa/2)$.

\begin{proposition}\label{prop:pc}
The function $\varphi_{\beta,m}(s)$ satisfies
\begin{equation}\label{def:varphi}
\varphi_{\beta,m}(s)=-\frac{1}{\sqrt{\pi}}\left( \frac{2}{\pi  }\right)^{s-\kappa/2} \frac{(s-1+\kappa/2)\Gamma(s+1/2-\kappa/2)}{(2s-1)\Gamma(s+1-\kappa/2)}
C(\beta,-m,s-\kappa/2).
\end{equation}
\end{proposition}

\begin{proof}
By Proposition \ref{FourierE} we have
\begin{equation}\label{eq:C}
C(\gamma,n,s)=
\frac{2^\kappa \pi^{s+\kappa}  n^{s+\kappa-1}}{\Gamma(s+\kappa)}\sum_{c\in\Z-\{0\}}|c|^{1-\kappa-2s} H_c^*(0,0,\gamma,n),
\end{equation}
if $n>0$. Comparing \eqref{def:b} and \eqref{eq:C} we find 
\[
b(0,0,s)=-C(\beta,-m,s-\kappa/2)\frac{s-1+\kappa/2}{(2s-1)\Gamma(s+1-\kappa/2)}.
\]
If we insert this into \eqref{eq:varphi}, we obtain the assertion.
\end{proof}

\subsection{Green's functions and meromorphic modular forms}

Let $\calC'\subset V_\C-\{0\}$ denote the cone over $\calK'\subset P(V_\C)$. 
Let $\Gamma\leq\Gamma(L)$ be a subgroup of finite index, $\chi$ a character\footnote{Throughout all characters are unitary.} of $\Gamma$, and $k$ an integer. By a meromorphic (holomorphic) modular form $F$ of weight $k$ with respect to $\Gamma$ and $\chi$ we mean a meromorphic (holomorphic) function on $\calC'$ satisfying
\begin{enumerate} 
\item $F$ is $\Gamma$-invariant, i.e., $F(\sigma Z)=\chi(\sigma) F(Z)$ for all $\sigma\in \Gamma$,
\item $F$ is homogeneous of degree $-k$, i.e., $F(a Z)= a^{-k} F(Z)$ for all $a\in \C\setminus\{0\}$.
\end{enumerate}
Our assumptions on $L$ from the beginning of section \ref{sect3} ensure that $F$ is automatically meromorphic (holomorphic) at the boundary of $X_L$ because of the Koecher principle. 

Observe that for $Z=X+iY\in \calC'$ we have
\[
(Z,\bar Z)= X^2+Y^2=2Y^2>0.
\]
If $F$ is a modular form of weight $k$ for the group $\Gamma$ with some character, then the function 
\begin{equation}\label{eq:petnorm}
\| F(Z)\| = |F(Z)|\cdot q(4\pi Y)^{k/2}
\end{equation}
on $\calC'$ is $\Gamma$-invariant and homogeneous of degree $0$. It is called the {\em Petersson norm} of $F$. 
It defines a Hermitean metric on the line bundle $\calL_k$ of modular forms of weight $k$.
The multiplicative normalization of the Petersson norm is adapted to the
normalization of Faltings heights in
arithmetic intersection theory. In particular it is chosen to be 
consistent with the normalization used in \cite{BrKue}, \cite{Kue}. 

The logarithm of the Petersson norm of any modular form, whose divisor consists of Heegner divisors, must be a linear combination of the Green's functions $\Phi_{\beta,m}(Z)$:

\begin{theorem}\label{thm:mf1}
Let $F$ be any meromorphic modular form of weight $k\in \Z$ for the group $\Gamma(L)$ with some character, and assume that
\begin{equation}\label{eq:divbor}
\div(F)=\frac{1}{2}\sum_{\beta\in L'/L}\sum_{\substack{m\in \Z+q(\beta)\\m<0}} a(\beta,m) H(\beta,m)
\end{equation}
is a linear combination of Heegner divisors.

i) The Petersson norm of $F$ is given by
\[
\log\|F(Z)\| = A -\frac{1}{8}\sum_{\beta\in L'/L}\sum_{\substack{m\in \Z+q(\beta)\\m<0}} a(\beta,m) \Phi_{\beta,m}(Z),
\]
where $A$ denotes an additive constant.

ii) The weight of $F$ is given by the coefficients of the Eisenstein series $E_0(\tau,0)$:
\[
k=  -\frac{1}{4} \sum_{\beta\in L'/L}\sum_{\substack{m\in \Z+q(\beta)\\m<0}}
a(\beta,m) C(\beta,-m,0).
\]
\end{theorem}

\begin{proof}
The first assertion is Theorem 4.23 in \cite{Br1}, and the second immediately follows from Corollary 4.24 there.
\end{proof}

Recall that a large class of modular forms for the group $\Gamma(L)$ with divisors consisting of Heegner divisors is given by {\em Borcherds products}  
(cf.~\cite{Bo2} Theorem 13.3,  \cite{Br1} Theorem 3.22). These are constructed as multiplicative liftings of $\C[L'/L]$-valued nearly holomorphic modular forms of weight $1-p/2$ for the group $\Mp_2(\Z)$ and the Weil representation $\rho_L$.
In fact, if the lattice $L$ splits two orthogonal hyperbolic planes over $\Z$, then every modular form with divisor as in \eqref{eq:divbor} is a constant multiple of a Borcherds product.

For our puroses it is important to notice that Borcherds product have nice arithmetic properties.
Their construction as infinte products, together with a result of McGraw \cite{MG}, implies that some integral power $\Psi^a$ of any Borcherds product $\Psi$ is a modular form,
whose Fourier coefficients at the cusp $z$ belong to $\Z[\zeta_N]$, where $\zeta_N$ denotes a primitive $N$-th root of unity. (The greatest common divisor of the Fourier coefficients of $\Psi^a$ at the cusp $z$ can be explicitly computed in terms of the coefficients of the Eisenstein series $E_\beta(\tau,0)$ with $\beta\neq 0$. It is equal to $1$, if $N=1$.) 
In particular, Borcherds products should define sections of the line bundle of modular forms over $\Z[\zeta_M]$ for some positive integer $M$, whenever there is a regular model of $X_L$ over $\Z[\zeta_M]$.

Finally, notice that there is a slight difference in the multiplicative normalizations of Borcherds products in \cite{Bo2} and \cite{Br1} (which only occurs if $N>1$). For our purposes here it is more convenient to work with Borcherds' normalization. To adapt \cite{Br1} to this normalization, one has to multiply
the generalized Borcherds product $\Psi_{\beta,m}(Z)$ in Definition 3.14 by the factor $\prod_{\ell=1}^{N-1} (1-e(\ell/N))^{b(\ell z/N,0)}$.
Then (3.40) still holds, if we replace the constant $C_{\beta,m}$ (see \eqref{eq:Cb}) by $\varphi'_{\beta,m}(\kappa/2)-2\varphi_{\beta,m}(\kappa/2)\Gamma'(1)$.

For questions regarding the arithmetic of Borcherds products, it is convenient to renormalize the Green's function $\Phi_{\beta,m}(Z)$ as follows.

\begin{definition}\label{def:G}
If $\beta\in L'/L$, and $m\in \Z+q(\beta)$ with $m<0$, then we define
\[
G_{\beta,m}(Z)=-\frac{1}{4}\big(\Phi_{\beta,m}(Z)-L_{\beta,m}\big),
\]
where 
\begin{equation}\label{def:L}
L_{\beta,m}=\varphi_{\beta,m}(\kappa/2)\left(\frac{\varphi'_{\beta,m}(\kappa/2)}{\varphi_{\beta,m}(\kappa/2)}-2\Gamma'(1)+\log(8 \pi^2)\right).
\end{equation}
\end{definition}

Proposition \ref{prop:aux} implies the identity
\[
G_{\beta,m}(Z)=-\frac{1}{4}\lim_{s\to \kappa/2}
\big(\Phi_{\beta,m}(Z,s)-2\zeta(2s+1-\kappa)(8\pi^2)^{s-\kappa/2}\varphi_{\beta,m}(s)\big).
\]

\begin{proposition}\label{thm:mf2}
If $F$ is a Borcherds product in the sense of \cite{Bo2} Theorem 13.3 with divisor as in \eqref{eq:divbor}, then 
\[
\log\|F(Z)\|^2 = \sum_{\beta\in L'/L}\sum_{\substack{m\in \Z+q(\beta)\\m<0}} a(\beta,m) G_{\beta,m}(Z).
\]
\end{proposition}

\begin{proof}
This immediately follows from Theorem 3.22 in \cite{Br1}.
\end{proof}

Proposition \ref{thm:mf2} can be viewed as
  a generalization of the well known Kronecker limit formula, 
which relates the special value of an eigenfunction of the Laplacian to 
the logarithm of the Petersson metric of an arithmetic modular form. 
To be more 
precise, if 
\[
\calE(\tau,s) = \frac{1}{2}\sum_{\gamma \in \Gamma_\infty
  \setminus \Gamma_1} (\Im (\gamma\tau))^s 
\]
is the real analytic Eisenstein series for $\Gamma_1=\Sl_2(\Z)$, then  the
logarithm of the Petersson norm 
of the Delta function is given by 
\begin{align}\label{eq:krolimdelta}
\log\left( |\Delta(\tau)|^2(4 \pi y)^{12} \right) =
- 4 \pi \lim_{s \to 1} \left( \calE(\tau,s) - \frac{\Gamma(1/2)
    \Gamma(s-1/2) \zeta(2s-1)}{\Gamma(s) \zeta(2s)}\right) + 12 \log(4 \pi).
\end{align}

We expect that the Green's function $G_{\beta,m}(Z)$ is a Green's
function  for the divisor 
$\frac{1}{2} H(\beta,m)$  in  the sense
of  \cite{BKK} which has a 
``good'' additive normalisation for questions in arithmetic
intersection theory (compare this with \cite{BrKue}).

\subsection{Integrals over Green's functions}\label{sect3.2}

The differential form  $\Omega=-dd^c\log(Z,\bar Z )$ on $\Gr(V)$ is $\Orth(V)$-invariant and positive. It is the associated $(1,1)$-form of an invariant K\"ahler metric $g$ on $\Gr(V)$. The corresponding invariant volume form  is given by $\frac{1}{p!}\Omega^p$.
We let
\[
B=\int\limits_{X_L}\Omega^p
\]
be the volume of $X_L$ and define the degree of 
a divisor $D$ on $X_L$ by
\[
\deg(D)=\int\limits_{D} \Omega^{p-1}.
\]

According to \eqref{eq:petnorm}, the function $|Y|^{-k}$ defines a Hermitean metric on the line bundle $\calL_k$ of modular forms of weight $k$ for $\Gamma(L)$, called the Petersson metric. The corresponding first Chern form of $\calL_k$ is given by
\[
c_1(\calL_k)=-k dd^c\log Y^2 = k\Omega.
\]

By means of results of Oda and Tsuzuki (based on an unfolding argument and the properties of spherical functions) it is possible to compute the integral of $\Phi_{\beta,m}(v)$ over $X_L$.

\begin{theorem}\label{thm:kron}
i) If $\Re(s)>\kappa/2$, then $\Phi_{\beta,m}(v,s)\in L^1(X_L)$. Moreover, if $f\in L^\infty(X_L)$ is a smooth bounded eigenfunction of the Laplacian $\Delta$ with eigenvalue $\lambda_f$, then 
\begin{equation}\label{eq:greenint}
\int\limits_{X_L} \Phi_{\beta,m}(v,s)f(v)\Omega^p
=
\frac{A(s)}{\lambda_f-\lambda_\Phi} 
\int\limits_{H(\beta,m)} f(v)\Omega^{p-1}.
\end{equation}
Here $A(s)$ denotes the normalizing factor
\begin{align}
\label{def:as}
A(s)= - \frac{p}{2\Gamma(s+1-\kappa/2)}
\end{align}
and
$\lambda_\Phi=\frac{1}{2}(s-\kappa/2)(s-1+\kappa/2)$ the eigenvalue of $\Phi_{\beta,m}(v,s)$. 

ii) The regularized Green's function $\Phi_{\beta,m}(v)$ belongs to $L^1(X_L)$ and
\begin{equation}\label{eq:kron}
\int\limits_{X_L} \Phi_{\beta,m}(v)\Omega^p=
\lim_{s\to\kappa/2}\int\limits_{X_L}\left(\Phi_{\beta,m}(v,s)-\frac{\varphi_{\beta,m}(\kappa/2)}{s-\kappa/2}\right)\Omega^p.
\end{equation}

\end{theorem}

\begin{proof}
The Green's function $\Phi_{\beta,m}(v,s)$ is a Poincar\'e series built out of the function
\begin{align}
\nonumber
g_s(v,\lambda)= 2\frac{\Gamma(s-1+\kappa/2)}{\Gamma(2s)} &\left(\frac{m}{m-q(\lambda_v)}\right)^{s-1+\kappa/2}\\
\label{def:g}
&\phantom{=}
\times F\left(s-1+\kappa/2,\, s+1-\kappa/2, \, 2s; \, \frac{m}{m-q(\lambda_v)}\right),
\end{align}
where $v\in \Gr(V)$ and $\lambda\in V$.
We now compare this function with the secondary spherical function $\phi^{(2)}_s(x)$ in \cite{OT}. As far as possible we use the notation of \cite{OT}.
We write $G_\R$ for the connected component of $\Orth(V)$ and fix a vector $\lambda^0$ in $\beta+L$ with $q(\lambda)=m$.  
There exists a basis $e_1,\dots,e_{p+2}$ of $V$ such that the quadratic form has the Gram matrix\footnote{Unspecified matrix entries are assumed to be $0$.} $2S$, where
\[
S=\begin{pmatrix} 1 &  &  &      &  \\
                    & 1&  &      &  \\
                    &  &-1&      &  \\
                    &  &  &\ddots&  \\
                    &  &  &      &-1\\
                  \end{pmatrix},
\]                  
and such that $\lambda^0=\sqrt{|m|}e_{p+2}$. 
We write $v_0$ for the element of $\Gr(V)$ given by the span of $ e_1,e_2$.
The stabilizer $K_\R$ in $G_\R$ of $v_0$ is a maximal compact subgroup.
Let $H_\R$ be the connected component of the stabilizer in $G_\R$ of $\lambda^0$. Clearly $ H_\R\cong \Orth^0(2,p-1)$.
We write $\Gamma$ for  the arithmetic subgroup of $G_\R$ given by the discriminant kernel $\Gamma(L)$ and put $\Gamma_H=\Gamma\cap H_\R$. We consider the function
\[
\calG_s(x,\lambda^0) = \sum_{\gamma\in \Gamma_H\bs\Gamma}  g_s(\gamma x v_0,\lambda^0)  
\qquad (x\in G_\R)
\]
and show that (up to a constant factor depending only on $s$) $\calG_{(s+1)/2}(x,\lambda^0)$ is equal to the function $G_s(x)$ defined in \cite{OT} \S3.1.
Since $\Phi_{\beta,m}(v,s)$ is a finite linear combination of the functions $\calG_s(x,\lambda^0)$, we will be able to apply the results of loc.~cit.

We begin by recalling some properties of the Lie algebra
$\frakg=\{ g\in M_{p+2}(\R);\; g^T S+Sg=0\}$ of $G_\R$ following \cite{OT} \S1.3.
The assignment $g\mapsto SgS$ defines a Cartan involution $\theta$ on $G_\R$ with fixed point subgroup $K_\R$.
If $J$ denotes the matrix
\[
J=\begin{pmatrix} 1 &  &      &  &  \\
                    & 1&      &  &  \\
                    &  &\ddots&  &  \\
                    &  &      & 1&  \\
                    &  &      &  &-1\\
                  \end{pmatrix},
\]
then $g\mapsto JgJ$ defines another involution $\sigma$ on $G_\R$ with fixed point subgroup $H_\R$. These involutions commute. They induce eigenspace decompositions $\frakg=\frakk\oplus \frakp$ (Cartan decomposition corresponding to $\theta$) and  $\frakg=\frakh\oplus \frakq$ (decomposition corresponding to $\sigma$) of $\frakg$. Here $\frakk$ denotes the Lie algebra of $K_\R$ and $\frakh$ the Lie algebra of $H_\R$.
We obtain the direct sum decomposition 
\[
\frakg=(\frakk\cap\frakh)\oplus (\frakk\cap \frakq)\oplus(\frakp\cap\frakh)\oplus(\frakp\cap\frakq).
\]
The subspace $\frakp\cap\frakq$ is given by matrices in $M_{p+2}(\R)$ of the form
\[
\begin{pmatrix}   &   &   &      &g_1     \\
                  & 0 &   &      &g_2     \\
                  &   &   &      & 0      \\
                  &   &   &      & \vdots \\
               g_1&g_2& 0 &\cdots& 0      \\
                  \end{pmatrix}
\]
with $g_1,g_2\in \R$.
Let 
\[
Y_0=
\begin{pmatrix}   &    &      & 1     \\
                  & 0  &      & 0     \\
                  &    &      & \vdots \\
                1 & 0  &\cdots& 0      \\
                  \end{pmatrix}
\]
and put $\fraka_{\frakp,\frakq}=\R Y_0$. Then $\fraka_{\frakp,\frakq}$ is a maximal Abelian subspace of $\frakp\cap\frakq$, and it is easily checked that we may take $Y_0$ for the $Y_0$ in \cite{OT}.

We now evaluate our function $g_s(x v_0,\lambda^0)$ (where $x\in G_\R$) on the $1$-parameter subgroup $\{ \exp(tY_0);\; t\in \R\}$ of $G_\R$.
To this end we have to determine the norm of
$(\exp(-tY_0)\lambda^0)_{v_0}$. Using the fact that 
\[
(Y_0^k e_{p+2})_{v_0}=
\begin{cases}
0,&2\mid k,\\
e_1,&2\nmid k,
\end{cases}
\]
we find that 
$(\exp(tY_0)e_{p+2})_{v_0}=e_1 \sinh(t)$. Hence 
\[
q\big((\exp(-tY_0)\lambda^0)_{v_0}\big)= m \sinh^2(t)
\]
and consequently
\begin{align*}
g_s\big(\exp(tY_0)v_0,\lambda^0\big)=
2&\frac{\Gamma(s-1+\kappa/2)}{\Gamma(2s)}\left(\cosh(t)\right)^{2-2s-\kappa}\\
&\times
F\left(s-1+\kappa/2,\, s+1-\kappa/2, \, 2s; \, \frac{1}{\cosh^2(t)}\right).
\end{align*}
On the other hand, by (2.5.3) of \cite{OT} (and noting $\rho_0=p/2$) we have
\begin{align*}
\phi^{(2)}_{2s-1}(\exp(tY_0))=
-&\frac{\Gamma(s-1+\kappa/2)\Gamma(s+1-\kappa/2)}{2\Gamma(2s)}
\left(\cosh(t)\right)^{2-2s-\kappa}\\
&\times
F\left(s-1+\kappa/2,\, s+1-\kappa/2, \, 2s; \, \frac{1}{\cosh^2(t)}\right).
\end{align*}

By virtue of Proposition 2.4.2. and \S2.5 of \cite{OT} we may infer
that
\begin{align}\label{compker}
g_s(x v_0,\lambda^0)&=A_0(s)\phi^{(2)}_{2s-1}(x),\\
\nonumber
A_0(s)&= -\frac{4}{\Gamma(s+1-\kappa/2)}.
\end{align}
Thus $\calG_s(x,\lambda^0) = A_0(s) G_{2s-1}(x)$, where $G_s(x)$ denotes the Green's function defined in \cite{OT} \S3.1 with respect to $H_\R$ and $\Gamma$.

Now Proposition 3.1.1 of \cite{OT} tells us that $\Phi_{\beta,m}(v,s)\in L^1(X_L)$ for $\Re(s)>\kappa/2$. 
To prove that $\Phi_{\beta,m}(v)\in L^1(X_L)$ and \eqref{eq:kron} one can use the explicit Fourier expansions with respect to the various parabolic subgroups (see \cite{Br1} Theorems 2.15 and 3.9) to derive estimates on Siegel domains following the argument of \cite{Br1} Chapter 4.2. 
Alternatively, one can use the fact that $\Phi_{\beta,m}(v,s)$ can be written as a regularized theta lifting together with a variant of the convergence proof of Kudla \cite{Ku2}. In fact, by means of the classical Weil estimate for the Siegel theta function and identity (2.20) in \cite{Br1} the claim can reduced to Propositions 3.2 and 3.4 of \cite{Ku2}.

We are left with proving \eqref{eq:greenint}.
Up to a constant factor the formula follows immediately from Proposition 3.3.1 in \cite{OT}. We have
\begin{equation}\label{eq:prel}
\frac{1}{B}\int\limits_{X_L} \Phi_{\beta,m}(v,s)f(v)\Omega^n
=\alpha
\frac{A_0(s)}{\lambda_f-\lambda_\Phi} 
\int\limits_{H(\beta,m)} f(v)\Omega^{p-1},
\end{equation}
with $\alpha \in \C$ coming from the different  normalizations of the Laplacian and the invariant measures on $G_\R$ and  $H_\R$.
We may use the fact that $\Phi_{\beta,m}(v,s)$ is an eigenfunction of $\Delta$ with eigenvalue $\lambda_\Phi$, the corresponding statement in \cite{OT} (see Proposition 2.4.2b), and \eqref{compker} to infer that the Laplacian in \cite{OT} equals $8$ times our $\Delta$.
Moreover, from identity (1.3.2) in \cite{OT} it can be concluded that
\[
\frac{\mu_{H_\R, OT}}{\vol_{OT}(X_L)}= p \cdot \frac{\mu_{H_\R, \Omega^{p-1}}}{\vol_{\Omega^p}(X_L)}.
\]
Here $\mu_{H_\R, OT}$ (respectively $\vol_{OT}(X_L)$) denotes the Haar measure on $H_\R$ (respectively the volume of $X_L$) normalized as in \cite{OT}, and $\mu_{H_\R, \Omega^{p-1}}$, $B=\vol_{\Omega^p}(X_L)$ the corresponding quantities in our normalzation. 
We find that $\alpha=p/8$.
If we insert $\alpha$ into \eqref{eq:prel} we obtain the assertion with $A(s)=\alpha A_0(s)$.
\end{proof}

M. Tsuzuki kindly
informs the authors that the meromorphic continuation of $\Phi_{\beta,m}(v,s)$ to a neighborhood of $s=\kappa/2$ and integrability results as in (ii) could also be obtained by means of the spectral theoretic methods developed in [OT].

\begin{proposition}\label{thm:deg}
i) The Eisenstein series $E_0(\tau,0)$ of weight $\kappa=1+p/2$ encodes the degrees of the Heegner divisors.
More precisely we have for $\kappa>2$:
\[
E_0(\tau,0)=2\frake_0-\frac{2}{B}\sum_{\gamma\in L'/L}
\sum_{\substack{n\in \Z-q(\gamma)\\n>0}}\deg(H(\gamma,-n))\frake_\gamma(n\tau).
\]
If $\kappa=2$, then the same identity holds, except for the fact that in the constant coefficient on the right hand side an additional term $\fraka y^{-1}$ might occur, where $\fraka\in \C[L'/L]$ is invariant under $\rho_L$.

ii) Using the notation of Theorem \ref{thm:eismain}, $\deg(H(\beta,-n))/B$ is equal to
\[
\begin{cases}
\displaystyle
(-1)^{\kappa/2-\delta/2} 2n^{\kappa-1} 
\frac{ |D_0|^{\kappa-1/2} }{\sqrt{|L'/L|}}
\cdot \frac{\sigma_{\gamma,n}(\kappa)}{L(\chi_{D_0},1-\kappa)},&\text{if $r$ is even,}\\[3ex]  
\displaystyle
(-1)^{\kappa/2-1/4+\delta/2}\frac{2^{2\kappa-3/2}n^{\kappa-1}|D_0|^{1-\kappa}}{\sqrt{|L'/L|}}
\cdot\frac{L(\chi_{D_0},3/2-\kappa)}{\zeta(2-2\kappa)} \sigma_{\gamma,n}(\kappa),&\text{if $r$ is odd.}   
\end{cases}
\]
\end{proposition}

\begin{proof}
If we apply \eqref{eq:greenint} of Theorem \ref{thm:kron} to the constant function $1$, we find
\begin{align*}
\int\limits_{X_L} \Phi_{\beta,m}(v,s)\Omega^p
&=-
\frac{A(s)}{\lambda_\Phi} 
\int\limits_{H(\beta,m)} \Omega^{p-1}\\
&=\frac{p}{(s-\kappa/2)(s-1+\kappa/2)}\cdot \frac{\deg(H(\beta,m))}{\Gamma(s+1-\kappa/2)}.
\end{align*} 
At $s=\kappa/2$ the right hand side has the Laurent expansion
\[
2\deg(H(\beta,m))(s-\kappa/2)^{-1} + O(1).
\]
But now it follows from \eqref{eq:kron} that the leading term is equal to $B\varphi_{\beta,m}(\kappa/2)(s-\kappa/2)^{-1}$.
Thus 
\[
\varphi_{\beta,m}(\kappa/2)=\frac{2}{B}\deg(H(\beta,m)).
\]
On the other hand, by \eqref{def:varphi} we know that $\varphi_{\beta,m}(\kappa/2)=-C(\beta,-m,0)$. This implies the first assertion. 
%

To prove the second statement we use the first assertion and Theorem \ref{thm:eismain}.  We get
\[
\frac{\deg(H(\beta,-n))}{B}=\begin{cases}\displaystyle
-  2n^{\kappa-1} 
\cos\big(\pi(-\tfrac{\delta}{2}+\tfrac{p-2}{4})\big)
\frac{|D_0|^{\kappa-1/2} }{\sqrt{|L'/L|}}
\cdot \frac{\sigma_{\gamma,n}(\kappa)}{L(\chi_{D_0},1-\kappa)},& \text{$r$ even,}\\[3ex]
\displaystyle
-\frac{2^{2\kappa-3/2}n^{\kappa-1}|D_0|^{1-\kappa}\sin(\pi\kappa)}{\sqrt{|L'/L|}\cos\big(\pi(-\tfrac{\delta}{2}+\tfrac{p-3}{4})\big)}
\cdot\frac{L(\chi_{D_0},3/2-\kappa)}{\zeta(2-2\kappa)} \sigma_{\gamma,n}(\kappa)
,& \text{$r$ odd.}
\end{cases}
\]
Recall that $\delta\in \{0,1\}$ is uniquely determined by $p$ modulo $4$ (see table on page \pageref{table}).
In the $r$ even case it is easily checked that 
\[
\cos\big(\pi(-\tfrac{\delta}{2}+\tfrac{p-2}{4})\big)=-(-1)^{\kappa/2-\delta/2}.
\]
Moreover, in the $r$ odd case we have
\[
\frac{\sin(\pi\kappa)}{\cos\big(\pi(-\tfrac{\delta}{2}+\tfrac{p-3}{4})\big)}
=-(-1)^{\kappa/2-1/4+\delta/2}.
\]
Putting this into the above equation, we obtain the assertion.
\end{proof}

To compute the coefficients $C(\gamma,n,0)$ and thereby the degrees of the $H(\beta,m)$ explicitly, it remains to determine the generalized divisor sums $\sigma_{\gamma,n}(\kappa)$. This can easily be carried out on a computer for any given lattice $L$.
A computer program can be downloaded from the first authors home-page.

\begin{theorem}\label{thm:maass}
Let $G_{\beta,m}(v)$ be the Green's function of Definition \ref{def:G}.  
If $h\in L^\infty(X_L)$ is a smooth bounded eigenfunction of the Laplacian $\Delta$ with eigenvalue $\lambda_h\neq 0$, then 
\[
\int\limits_{X_L} G_{\beta,m}(v) h(v)\Omega^p 
=\frac{p}{8\lambda_h}
\int\limits_{H(\beta,m)} h(v)\Omega^{p-1}.
\]
In particular, if $F$ is a Borcherds product in the sense of \cite{Bo2} Theorem 13.3, then 
\[
\int\limits_{X_L} \log\|F(v)\|^2 h(v)\Omega^p 
=\frac{p}{8\lambda_h}
\int\limits_{\dv(F)} h(v)\Omega^{p-1}.
\]
\end{theorem} 

\begin{proof}
The symmetry of the Laplacian and the assumption $\lambda_h\neq 0$ imply that $\int_{X_L}h(v)\Omega^p=0$. Therefore, inserting the definition of $G_{\beta,m}(v)$ we find
\begin{align*}
\int\limits_{X_L} G_{\beta,m}(v) h(v)\Omega^p 
&=-\frac{1}{4}\int\limits_{X_L} (\Phi_{\beta,m}(v)-L_{\beta,m}) 
h(v)\Omega^p \\
&=-\frac{1}{4}\int\limits_{X_L} \Phi_{\beta,m}(v)
h(v)\Omega^p .
\end{align*}
As in \eqref{eq:kron} this can be written in the form
\begin{align*}
&-\frac{1}{4}\lim_{s\to\kappa/2} \int\limits_{X_L} \left(
\Phi_{\beta,m}(v,s)-\frac{\varphi_{\beta,m}(\kappa/2)}{s-\kappa/2}\right)
h(v)\Omega^p\\
&=-\frac{1}{4}\lim_{s\to\kappa/2} \int\limits_{X_L}
\Phi_{\beta,m}(v,s)h(v)\Omega^p.
\end{align*}
Hence, by means of Theorem \ref{thm:kron} we obtain
\begin{align*}
\int\limits_{X_L} G_{\beta,m}(v) h(v)\Omega^p 
&=\frac{p}{8\lambda_h}
\int\limits_{H(\beta,m)} h(v)\Omega^{p-1}.
\end{align*}
The second assertion follows from Proposition \ref{thm:mf2}.
\end{proof}

If $h\in L^\infty(X_L)$ is a bounded eigenfunction of the Laplacian $\Delta$ with eigenvalue $\lambda_h=0$, then $h$ has to be constant.

\begin{theorem}\label{thm:int}
Let $G_{\beta,m}(v)$ be the Green's function of Definition \ref{def:G}.  
Then
\[
\frac{1}{B}\int\limits_{X_L} G_{\beta,m}(v) \Omega^p 
=-\frac{C(\beta,-m,0)}{4}
\left(\frac{C'(\beta,-m,0)}{C(\beta,-m,0)}+\log(4\pi)-\Gamma'(1)\right).
\]
Here $C(\gamma,n,s)$ denotes the $(\gamma,n)$-th  coefficient of the Eisenstein series $E_0(\tau,s)$ of weight $\kappa=1+p/2$.
In particular, if $F$ is a Borcherds product of weight $k$ in the sense of \cite{Bo2} Theorem 13.3 with divisor as in \eqref{eq:divbor}, then 
\[
\frac{1}{B}\int\limits_{X_L} \log\|F(v)\|^2 \Omega^p 
=k\left(\log(4\pi)-\Gamma'(1)\right)
-\frac{1}{4} \sum_{\beta, m} a(\beta,m)C'(\beta,-m,0).
\]
\end{theorem} 

\begin{proof}
We use Theorem \ref{thm:kron} to compute the integral:
\begin{align*}
\frac{1}{B}\int\limits_{X_L} G_{\beta,m}(v) \Omega^p 
&=-\frac{1}{4B}\int\limits_{X_L}\left(\Phi_{\beta,m}(v)-L_{\beta,m}\right)
\Omega^p\\
&=-\frac{1}{4B}\lim_{s\to\kappa/2}\int\limits_{X_L}\left(\Phi_{\beta,m}(v,s)-\frac{\varphi_{\beta,m}(\kappa/2)}{s-\kappa/2}\right)\Omega^p
+\frac{L_{\beta,m}}{4} \\
&=\lim_{s\to\kappa/2}\left( \frac{A(s)}{4B\lambda_\Phi}\deg(H(\beta,m))+\frac{\varphi_{\beta,m}(\kappa/2)}{4(s-\kappa/2)}\right)+\frac{L_{\beta,m}}{4}.
\end{align*}
If we insert the Laurent series expansions of $\lambda_\Phi=\frac{1}{2}(s-\kappa/2)(s-1+\kappa/2)$ and 
\begin{align*}
A(s)&=-\frac{p}{2}+\frac{p}{2}\Gamma'(1)(s-\kappa/2) + O((s-\kappa/2)^2)
\end{align*}
at $s=\kappa/2$, and use that $2\deg(H(\beta,m))=B\varphi_{\beta,m}(\kappa/2)$,
we obtain
\begin{align*}
\frac{1}{B}\int\limits_{X_L} G_{\beta,m}(v) \Omega^p 
&=\frac{\varphi_{\beta,m}(\kappa/2)}{4}\left(\Gamma'(1)+\frac{1}{\kappa-1}\right)+\frac{L_{\beta,m}}{4}.
\end{align*}
Inserting the definition \eqref{def:L} of $L_{\beta,m}$,
we get 
\begin{equation}\label{eq:mo}
\frac{1}{B}\int\limits_{X_L} G_{\beta,m}(v) \Omega^p 
=\frac{\varphi_{\beta,m}(\kappa/2)}{4}\left(\frac{\varphi_{\beta,m}'(\kappa/2)}{\varphi_{\beta,m}(\kappa/2)}+\frac{1}{\kappa-1}-\Gamma'(1)+\log(8\pi^2)\right).
\end{equation}
It follows from \eqref{def:varphi} that
\begin{align*}
\varphi_{\beta,m}(\kappa/2)&=-C(\beta,-m,0),\\
\frac{\varphi'_{\beta,m}(\kappa/2)}{\varphi_{\beta,m}(\kappa/2)}
&=\frac{C'(\beta,-m,0)}{C(\beta,-m,0)}-\log(2\pi)-\frac{1}{\kappa-1}.
\end{align*}
Putting this into \eqref{eq:mo}, we obtain the first assertion, the second follows from Proposition \ref{thm:mf2} and Theorem \ref{thm:mf2} (ii).
\end{proof}

We are now ready to prove the main result of this section.
Recall that $G_{\beta,m}(Z)$ is an automorphic Green's function for $\frac{1}{2}H(\beta,m)$. 

\begin{theorem}\label{thm:main}
i) If $G_{\beta,m}(v)$ denotes the Green's function of Definition \ref{def:G}, then  
\begin{equation}\label{eq:t}
\frac{2}{\deg H(\beta,m)}\int\limits_{X_L} G_{\beta,m}(v) \Omega^p 
=
\frac{C'(\beta,-m,0)}{C(\beta,-m,0)}+\log(4\pi)-\Gamma'(1).
\end{equation}
Here $C(\gamma,n,s)$ denotes the $(\gamma,n)$-th  coefficient of the Eisenstein series $E_0(\tau,s)$ of weight $\kappa=1+p/2$.

ii) Using the notation of Theorem \ref{thm:eismain}, this is equal to
\[
2\frac{L'(\chi_{D_0},1-\kappa)}{L(\chi_{D_0},1-\kappa)}+2\frac{\sigma_{\beta,-m}'(\kappa)}{\sigma_{\beta,-m}(\kappa)}+\log|mD_0^2|+\sum_{j=1}^{\kappa-1}\frac{1}{j},
\]
if $r$ even, and to
\[
4\frac{\zeta'(2-2\kappa)}{\zeta(2-2\kappa)}-2\frac{L'(\chi_{D_0},3/2-\kappa)}{L(\chi_{D_0},3/2-\kappa)}
+2\frac{\sigma_{\beta,-m}'(\kappa)}{\sigma_{\beta,-m}(\kappa)}+\log|4m/D_0^2|
+\sum_{j=1}^{\kappa-1/2}\frac{2}{2j-1},
\]
if $r$ is odd.
\end{theorem}

\begin{proof}
The first statement is an immediate consequence of Theorem \ref{thm:int} and  Proposition \ref{thm:deg}.

If we compute the logarithmic derivative of $C(\beta,-m,s)$ at $s=0$  by means of Theorem \ref{thm:eismain}, we find that it is equal to
\[
2\frac{L'(\chi_{D_0},1-\kappa)}{L(\chi_{D_0},1-\kappa)}+2\frac{\sigma_{\beta,-m}'(\kappa)}{\sigma_{\beta,-m}(\kappa)}-\log(4\pi)+\log|m|+\log(D_0^2)+\frac{\Gamma'(\kappa)}{\Gamma(\kappa)},
\]
if $r$ even, and to
\[
4\frac{\zeta'(2-2\kappa)}{\zeta(2-2\kappa)}-2\frac{L'(\chi_{D_0},3/2-\kappa)}{L(\chi_{D_0},3/2-\kappa)}
+2\frac{\sigma_{\beta,-m}'(\kappa)}{\sigma_{\beta,-m}(\kappa)}-\log(\pi/4)+\log|m|-\log(D_0^2)+\frac{\Gamma'(\kappa)}{\Gamma(\kappa)},
\]
if $r$ odd. Here we have used that the logarithmic derivatives of the $\sin$- and $\cos$-terms vanish at $s=0$. 
The functional equation of the $\Gamma$-function implies that
\[
\frac{\Gamma'(\kappa)}{\Gamma(\kappa)}=\begin{cases}
\Gamma'(1)+\sum_{j=1}^{\kappa-1}\frac{1}{j},&\text{if $r$ is even}\\
\Gamma'(1)-2\log(2)+\sum_{j=1}^{\kappa-1/2}\frac{1}{j-1/2},&\text{if $r$ is odd.}
\end{cases}
\]
Inserting this into the above formula and then into \eqref{eq:t}, we obtain the assertion.
\end{proof}

If we rescale the Petersson metric by a constant factor $a$, then, by Proposition \ref{thm:deg}, the right hand side of \eqref{eq:t} changes by the additive constant $\log(a)$. In that way one could absorb the constant $\log(4\pi)-\Gamma'(1)$ in \eqref{eq:t}. Moreover, one could compensate changes in the formulas resulting from rescaling $E_0(\tau,s)$ by a factor depending on $s$.

\section{Examples}
\label{sect4}

{\bf 1. }
We first consider the case of the Siegel modular group of genus $2$ and confirm the result of Kudla \cite{Ku2}.
Let $V=\R^5$ equipped with the quadratic form 
\[
q((x_1,\dots,x_5))=x_1 x_2+x_3 x_4 -x_5^2.
\]
Then $(V,q)$ is a quadratic space of signature $(2,3)$ and $\kappa=r/2=5/2$.
Furthermore, $L=\Z^5$ is an even lattice in $V$, whose discriminant group $L'/L$ is isomorphic to $\Z/2\Z$.
The space of modular forms of weight $\kappa$ with respect to the Weil representation $\rho_L^*$ and $\Mp_2(\Z)$ is isomorphic to the space of modular forms of weight $\kappa$ for the group $\Gamma_0(4)$ having a Fourier expansion of the form
$
f(\tau)=\sum_{n=0}^\infty a(n)e(n\tau)
$
with $a(n)=0$ if $n\not\equiv 0,1\pmod{4}$
(Kohnen's plus space \cite{Ko}). Thus the Eisenstein series $E_0(\tau,s)$ defined in \eqref{DefEb} is essentially the non-holomorphic Cohen Eisenstein series of weight $\kappa$. 

Let $D$ be a positive fundamental discriminant. We simply write $C(D,s)$ for the coefficient $C(\beta,D/4,s)$ of $E_0(\tau,s)$, where $\beta\in L'/L$ is taken to be $0+L$, if $D\equiv 0\pmod{4}$, and the non-zero element of $L'/L$, if $D\equiv 1\pmod{4}$.
We may compute $C(D,s)$ by means of Theorem \ref{thm:eismain}. Noting that $D_0=D$ and $\delta=0$,
we find
\begin{align*}
C(D,s)&=
2\pi^{-s}D^{-s} 
\frac{\sin(\pi(2s+\kappa))\Gamma(2s+\kappa)}{\cos(\pi s)\Gamma(s+\kappa)}
\cdot\frac{L(\chi_{D},3/2-2s-\kappa)}{\zeta(2-4s-2\kappa)} \sigma_{\beta,D/4}(2s+\kappa)\\
&=2\pi^{-s}D^{-s} 
\frac{\cos(2\pi s)\Gamma(2s+5/2)}{\cos(\pi s)\Gamma(s+5/2)}
\cdot\frac{L(\chi_{D},-1-2s)}{\zeta(-3-4s)} \sigma_{\beta,D/4}(2s+5/2),
\end{align*}  
where $\sigma_{\beta,D/4}(s)$ is given by \eqref{eq:defsigma}.
In view of Proposition \ref{prop:genus} it is easily verified that $\sigma_{\beta,D/4}(s)=1$. Here we need the assumption that $D$ be fundamental. For general discriminants $\sigma_{\beta,D/4}(s)$ can be computed using 
\cite{EZ} \S2 p.~21 or \cite{BK} Example 10. The result is $\sigma_{\beta,D/4}(s)=\sum_{d\mid f} \mu(d)\chi_{D_0}(d)d^{1/2-s}\sigma_{2-2s}(f/d)$, where $D=D_0 f^2$ as in section \ref{sect2}.

It is well known that in this case the discriminant kernel $\Gamma(L)$ is isomorphic to the Siegel modular group $\Symp_2(\Z)$ of genus $2$ (see for instance \cite{GN}). The quotient $X_L$ is isomorphic to the Siegel modular threefold $\Symp_2(\Z)\bs \H_2$, where $\H_2$ denotes the Siegel upper half plane of genus $2$.
The Heegner divisor $\frac{1}{2}H(\beta,-D/4)$ ($\beta$ being determined by $D$) can be identified with the Humbert surface $H(D)$ of discriminant $D$ (cf.~\cite{Ge} Chapter IX). It is an irreducible divisor on $\Symp_2(\Z)\bs \H_2$, which is birational to the symmetric Hilbert modular surface of discriminant $D$.

In the coordinates of $\H_2$ the K\"ahler form $\Omega$ is given by $-dd^c \log(\det(Y))$, where $Z\in \H_2$ and $Y=\Im(Z)$. It is well known that 
$B=\int_{\Symp_2(\Z)\bs \H_2}\Omega^3= -\zeta(-1)\zeta(-3)$.
According to Proposition \ref{thm:deg} we have
\[
\deg(H(D))=-\frac{B}{4}C(D,0)=-\frac{B}{2} 
\frac{L(\chi_{D},-1)}{\zeta(-3)}=\frac{1}{2}\zeta_K(-1),
\]
where $\zeta_K(s)$ denotes the Dedekind zeta function of $K=\Q(\sqrt{-1})$; 
a result which was first proved by van der Geer (see \cite{Ge} Chapter IX.2).

For any Humbert surface $H(D)$ there exists a Borcherds product $\Psi_D(Z)$ with divisor $H(D)$. It can be viewed as the Borcherds lifting (\cite{Bo2} Theorem 13.3) of the unique modular form of weight $-1/2$ for $\Gamma_0(4)$, which is holomorphic on $\H$ and whose Fourier expansion has the form 
\[
f(\tau)=q^{-D}+\sum_{\substack{n\geq0\\-n\equiv 0,1\;(4)}}^\infty a(n) e(n\tau).
\]
If we briefly write $G_D(Z)$ for the Green's function $G_{\beta,-D/4}(Z)$ of Definition \ref{def:G}, then by Proposition \ref{thm:mf2} we have 
\[
\log\|\Psi_D(Z)\|^2=\log \left( |\Psi_D(Z)|^2 (\det(4\pi Y))^{k}\right)=G_D(Z),
\]
where $k=-C(D,0)/4$ is the weight of $\Psi_D(Z)$.
By virtue of Theorem \ref{thm:main} we find that
\begin{align*}
&\frac{2}{\zeta_K(-1)}\int\limits_{\Symp_2(\Z)\bs \H_2}\log\|\Psi_D(Z)\|^2\,\Omega^3=
4\frac{\zeta'(-3)}{\zeta(-3)}-2\frac{L'(\chi_D,-1)}{L(\chi_D,-1)}+\frac{8}{3}
-\log(D).
\end{align*}

\bigskip

{\bf 2. } 
We now consider Hilbert modular surfaces using our results for the orthogonal group $\Orth(2,2)$. This case was also considered in \cite{BrKue}, but working entirely with the group $\Sl_2(\OK)$ and the Green's functions as defined in \cite{Br2}.  

Let $D>0$ be a positive fundamental discriminant and $K=\Q(\sqrt{D})$ the real quadratic field of discriminant $D$. For simplicity we assume that $D$ is a prime.
We write $a\mapsto a'$ for the conjugation, $\norm(a)=aa'$ for the norm, and $\OK$ for the ring of integers in $K$. As usual we view the corresponding Hilbert modular group $\Gamma_K=\Sl_2(\OK)$ as a subgroup of $\Sl_2(\R)\times \Sl_2(\R)$. We briefly recall some facts on the identification of $(\Sl_2(\R)\times\Sl_2(\R))/\{\pm (1,1)\} $ with the orthogonal group $\Orth^0(2,2)$ (for more details see \cite{Ge} Chapter V.4 and \cite{Bo3} Example 5.5).

Let $V$ be the vector space of real $2\times 2$ matrices equipped with the quadratic form given by $q(M)=-\det(M)$. Then $(V,q)$ is a real quadratic space of signature $(2,2)$ and $\kappa=r/2=2$.
We consider the even lattice $L\subset V$ of matrices $X=\kzxz{a}{\nu}{\nu'}{b}$ with $a,b\in \Z$ and $\nu\in \OK$. The dual lattice $L'$ 
is given by matrices $X$ as before, but with $\nu\in \frakd^{-1}$, the inverse of the different.
The group $\Gamma_K$ acts on $V$ by $X\mapsto MX {M'}^t$ for 
$M\in \Gamma_K$, the quadratic form and the lattice $L$ being preserved. In that way one gets an isomorphism $\Gamma_K/\{\pm 1\}\to \Gamma(L)$.
The Grassmannian $\Gr'(V)$ can be identified with the product $\H^2$ of two copies of $\H$. The action of $\Sl_2(\R)\times\Sl_2(\R)$ by fractional linear transformations corresponds to the linear action of $\Orth'(V)$ on $\Gr'(V)$.
Therefore modular forms for $\Gamma(L)$ on $\Gr'(V)$ can be identified with Hilbert modular forms on $\H^2$ for the group $\Gamma_K$. Moreover, the Heegner divisors on $X_L$ correspond to Hirzebruch-Zagier divisors on $\Gamma_K\bs \H^2$.
In the coordinates of $\H^2$ we have 
\[
\Omega=-dd^c\log(y_1 y_2)= \frac{1}{4 \pi} \left( \frac{dx_1   dy_2}{y_1^2}
+ \frac{dx_2   dy_2}{y_2^2} \right),
\] 
where $Z=(z_1,z_2)\in\H^2$, and $(y_1,y_2)=\Im(z_1,z_2)$. It is well known that $B=\int_{\Gamma_K\bs\H^2}\Omega^2=\zeta_K(-1)$, where $\zeta_K(s)$ denotes the Dedekind zeta function of $K$.

Let $\beta\in L'/L=\frakd^{-1}/\OK\cong\Z/D\Z$ and $M$ a positive integer such that $M/D\in \Z-q(\beta)$, that is $M\equiv -\norm(\beta)\pmod{D}$.
If we compute the coefficient $C(\beta,M/D,s)$ of the Eisenstein series $E_0(\tau,s)$ of weight $\kappa$ using Theorem \ref{thm:eismain} (noting $D_0=D$), we find that
\begin{align*}
C(\beta,M/D,s)&=2^{2-2s} \pi^{-s}  D^{s}M^{\kappa+s-1} 
\cos(\pi s)
\frac{\Gamma(2s+\kappa)}{\Gamma(s+\kappa)}
\cdot \frac{\sigma_{M}(2s+\kappa)}{L(\chi_{D},1-2s-\kappa)},
\end{align*}
where
\begin{equation}\label{eq:chil}
\sigma_M(s)=\prod_{p\mid 2MD} \frac{L^{(p)}_{\beta,M/D}(p^{1-r/2-s})}{1-\chi_D(p)p^{-s}}.
\end{equation}
To determine the $p$-polynomials in the finite Euler product $\sigma_M(s)$ more explicitly one can use Proposition  \ref{prop:genus}, which reduces the computation of representation numbers of the lattice $L$ modulo prime powers to the computation of such numbers for the smaller lattice $\OK$.
We see that
\[
L^{(p)}_{\beta,M/D}(p^{-r/2}X)
=\tilde N_{\beta,M/D}(p^{w_p})(p^{-1}X)^{w_p}+(1-X)\sum_{\lambda=0}^{w_p-1}\tilde N_{\beta,M/D}(p^\lambda)(p^{-1}X)^\lambda,
\]
where
\[
\tilde N_{\beta,M/D}(p^\lambda)=\#\{ x\in \OK/p^\lambda\OK; \quad \norm(x-\beta)+M/D\equiv 0\pmod{p^\lambda}\}.
\]
It is easily verified that
\[
\tilde N_{\beta,M/D}(p^\lambda)=
\begin{cases}\#\{ y\in \OK/p^\lambda\OK; \quad \norm(y)\equiv M\pmod{p^\lambda}\},&
p\neq D,\\
\frac{S(M)}{D}\cdot\#\{ y\in \OK/p^{\lambda+1} \OK; \quad \norm(y)\equiv M\pmod{p^{\lambda+1}}\},& p =D,
\end{cases}
\]
where $S(M)=1$, if $D$ divides $M$, and $S(M)=1/2$, if $D$ does not divide $M$.
If we define 
\[
N_b(M)=\# \{ x\in \OK/b\OK;\quad  \norm(x)\equiv M\pmod{b}\}
\]
as in \cite{Za2} p.~27, we find that $L^{(p)}_{\beta,M/D}(p^{-r/2}X)$ is equal to
\[
\begin{cases}
N_{p^{w_p}}(M)(p^{-1}X)^{w_p}+(1-X)\sum_{\lambda=0}^{w_p-1}N_{p^\lambda}(M)(p^{-1}X)^\lambda,&p\neq D,\\
\frac{S(M)}{D}\left(
N_{p^{w_p+1}}(M)(p^{-1}X)^{w_p}+(1-X)\sum_{\lambda=0}^{w_p-1}N_{p^{\lambda+1}}(M)(p^{-1}X)^\lambda\right),&
p= D.
\end{cases}
\]
The representation numbers $N_{p^\lambda}(M)$ are computed in Lemma 3 on p.~27 of \cite{Za2}. This can be used to determine  $L^{(p)}_{\beta,M/D}(p^{-r/2}X)$ more explicitly. Since the computation is somewhat lengthy (but trivial), we omit the details. If we write $M=M_0 D^\nu$ with $M_0$ coprime to $D$, we obtain
\[
L^{(p)}_{\beta,M/D}(p^{-r/2}X)=\begin{cases}
\big(1-\chi_D(p)p^{-1} X\big)\sum_{\lambda=0}^{v_p(M)} \big(\chi_D(p)X\big)^\lambda,& p\neq D,\\
S(M)(1+\chi_D(M_0)X^\nu),&p=D.
\end{cases}
\]
Inserting this into \eqref{eq:chil} we infer
\begin{align*}
\sigma_M(s)
&=S(M)(1+\chi_D(M_0)D^{(1-s)\nu})\prod_{p\mid 2MD} \sum_{\lambda=0}^{v_p(M)} \chi_D(p^\lambda)p^{(1-s)\lambda}\\
&=S(M)\sum_{d\mid M}\big( \chi_D(d)+\chi_D(M/d)\big) d^{1-s},
\end{align*}
and therefore
\begin{align*}
C(\beta,M/D,s)&=2^{2-2s} \pi^{-s}  D^{s}M^{s+1} 
\frac{\cos(\pi s)
\Gamma(2s+2)S(M)}{\Gamma(s+2)L(\chi_{D},-1-2s)}
\sum_{d\mid M}\big( \chi_D(d)+\chi_D(M/d)\big) d^{-1-2s}.
\end{align*}

The Hirzebruch Zagier divisor $T(M)$ on $\Gamma_K\bs \H^2$  (cf.~\cite{Ge} Chap.~V) can be identified with 
the divisor $\frac{1}{2}H(\beta,-M/D)$, if $D$ divides $M$, and with $H(\beta,-M/D)$, 
if $D$ does not divide $M$. The multiplicities of the irreducible components of $T(M)$ are equal to $1$.
Therefore we define the Green's function on $\H^2$ corresponding to $T(M)$ by 
$G_M(Z)=S(M)^{-1} G_{\beta,-M/D}(Z)$. 
It has a logarithmic singularity  along $T(M)$ of type $\log|f|^2$, where $f$ denotes a local holomorphic equation for $T(M)$. 
Now Proposition \ref{thm:deg} says that
\[
\deg(T(M))=-\frac{B}{4S(M)}C(\beta,M/D,0)=-\zeta(-1)
\sum_{d\mid M}\big( \chi_D(d)+\chi_D(M/d)\big) d,
\]
a result which was already proved by Hirzebruch-Zagier.
By means of Theorem \ref{thm:main} we obtain 
\begin{align} 
&\frac{1}{\deg(T(M))} \int\limits_{\Gamma_K\bs \H^2} G_M(Z) \Omega^2=2
\frac{L'(\chi_D,-1)}{L(\chi_D,-1)}+1+2\frac{\sigma_M'(2)}{\sigma_M(2)}+\log(DM).
\label{eq:hilex}
\end{align}
If $\Psi(Z)$ is a Borcherds product (see Theorem 9 in \cite{BB}), then its divisor is a linear combination of the $T(M)$ and the integral over $\log\|\Psi(Z)\|^2$ is given by the corresponding linear combination of the quantities on the right hand side of \eqref{eq:hilex}.

For instance, if $D=5$ and $K=\Q(\sqrt{D})$, the product $\Theta$ of the $10$ theta constants considered by Gundlach \cite{Gu} 
is a Hilbert modular form of weight $5$ for the group $\Gamma_K$ with divisor $T(1)$. The Fourier coefficients of $\Theta$ are integral and have greatest common divisor $64$.
Thus $2^{-6}\Theta$ is a Borcherds product and 
\[
\log\|2^{-6}\Theta(z_1,z_2)\|^2=\log\left( |2^{-6}\Theta(z_1,z_2)|^2(16\pi^2 y_1 y_2)^{5}\right)=G_1(z_1,z_2).
\]
Hence, inserting $\sigma_1(s)=1$ and $\deg(T(1))=-2\zeta(-1)$, we get
\[
\frac{-1}{2\zeta(-1)}\int\limits_{\Gamma_K\bs \H^2}\log\|2^{-6}\Theta(z_1,z_2)\|^2\Omega^2=
 2\frac{L'(\chi_D,-1)}{L(\chi_D,-1)}+1+\log(D).
\]

\bigskip

{\bf 3. } 
We finally consider the classical modular group $\Gamma_1=\Sl_2(\Z)$ including the example \eqref{eq:werbung} mentioned in the introduction. 
This could be done using the exceptional isomorphism relating $\Sl_2(\R)$ to the orthogonal group $\Orth(2,1)$ and by carefully extending our results (and some of \cite{Br1}) to $p=1$. However, because of the difficulties caused by possible singularities at the cusps and by convergence questions, we chose to give a direct proof based on Rohrlich's modular Jensen formula and the classical Kronecker limit formula.

Let $D$ be a negative fundamental discriminant and $K=\Q(\sqrt{D})$. We
 briefly recall some properties of  Heegner divisors. Every ideal class $[\fraka]$ of $K$ defines a unique point $[\rho_\fraka]$ on $\Gamma_1\setminus \H$ 
by associating with a fractional ideal $\fraka=\Z a+\Z b$ with oriented (i.e.~$\Im(b\bar a)>0$) $\Z$-basis $a,b$ the point $\rho_\fraka=b/a\in \H$.
The Heegner divisor $H(D)$ on $\Gamma_1\setminus \H$ consists of the sum of the $[\rho_\fraka]$, where $\fraka$ runs through all ideal classes of $K$, counted with multiplicity $2/w$, where $w$ is the number of units in $K$.
 The cardinality of $H(D)$ is equal to the class number $h$ of $K$, its degree is  $2 h/ w$.
We write $\Psi_D(\tau)$ for the unique modular form for $\Gamma_1$ whose divisor equals $H(D)$ and whose value at the cusp 
$\infty$ is given by  $\Psi_D(\infty)=1$.

\begin{theorem} 
Let $D$ and $H(D)$ be a as above. The
  degree of $H(D)$ is equal to $L(\chi_D,0)$ and 
\begin{align}\label{eq:beispiel1}
\frac{1}{L(\chi_D,0)}  \int\limits_{\Gamma_1\setminus \H}
 \log\left( |\Psi_D(\tau)|^2 (4 \pi
   y)^{k}\right) \frac{dxdy}{4 \pi y^2} = 2 \frac{\zeta'(-1)}{\zeta(-1)}+ 1
-  \frac{L'(\chi_D,0)}{L(\chi_D,0)}  -\frac{1}{2}   \log|D|.    
\end{align}
Here $k=24h/w$ denotes the weight of $\Psi_D$.
\end{theorem}

\begin{proof}
If $f \in M_k(\Gamma_1)$ is a modular of weight $k$ for $\Gamma_1$ with
$f(\infty)=1$, then 
 a beautiful formula of Rohrlich \cite{Ro} says: 
\begin{align*}
&\frac{3}{\pi} \int\limits_{\Gamma_1\setminus \H} \log\left( |f(\tau)|^2 (4 \pi
  y)^{k}\right) \frac{dxdy}{ y^2} = 
 k   \left(
 2  \frac{\zeta'(-1)}{\zeta(-1)}+1\right)
-  \sum_{\rho \in   \dv(f)}  m_\rho \log\left( |\Delta(\rho)|^2 (4 \pi
  \Im \rho )^{12}\right).
\end{align*}
In the sum we view $\dv(f)= \sum_\rho m_\rho\, [\rho]$ as a divisor on
$\Sl_2(\Z)\setminus \H$ such that $\sum m_\rho = k/12$. 
Recall the Kronecker limit formula \eqref{eq:krolimdelta}
and the identity
\begin{align} \label{eq:CM-keyformula} 
\sum_{\rho \in H(D)} \calE( \rho, s) = 
\frac{w}{2} \left|\frac{D}{4}\right|^{s/2} \frac{\zeta_{K}(s)}{  \zeta(2s)}, 
\end{align}
where $\zeta_{K}(s)= \zeta(s) L(\chi_D,s)$
denotes the Dedekind zeta function of $K$ (see \cite{GZ} p.~210). 
By means of the functional equations we obtain the Laurent expansions 
\begin{align*}
\frac{1}{\zeta(2s)} &= 
\frac{-1}{2\pi^2 \zeta(-1)}\left(
1+ \left(2\frac{\zeta'(-1)}{\zeta(-1)} +2  -2 \gamma- 2 \log (\pi)-2 \log(2)+
  \right) (s-1)\right)+\ldots,\\
\zeta_K(s)&=\frac{\pi L(\chi_D,0)}{\sqrt{D}} \left(
 (s-1)^{-1}
 -\frac{L'(\chi_D,0)}{L(\chi_D,0)}+2 \gamma+
 \log(2)-\log|D|+\log(\pi)
 \right)+\ldots.
\end{align*}
Hence the Eisenstein series $\calE(\tau,s)$ has the expansion
\begin{align*} 
\sum_{\rho \in H(D)} \calE( \rho, s)=& \frac{3 w L(\chi_D,0)}{2 \pi}
\left( (s-1)^{-1}+ 2 \frac{\zeta'(-1)}{\zeta(-1)}
  -\frac{L'(\chi_D,0)}{L(\chi_D,0)} +2 -
\frac{1}{2} \log|D| -\log(4 \pi) \right)\\
 &+O(s-1)
\end{align*}
at $s=1$. Finally we have
$$
 \frac{\Gamma(1/2)
    \Gamma(s-1/2) \zeta(2s-1)}{\Gamma(s) \zeta(2s)} =
\frac{3}{\pi}
\left( (s-1)^{-1}+   2 \frac{\zeta'(-1)}{\zeta(-1)} +2 - 2 \log (4
  \pi)
 \right)+O(s-1).
$$
By the Kronecker limit formula we derive the degree relation
\begin{align}\label{eq:classnumber}
 \frac{w  L(\chi_D,0)}{2} = h.
\end{align}
Moreover, we find
\begin{align*}
&\sum_{\rho \in H(D)}  \log\left( |\Delta(\rho)|^2(4 \pi \Im \rho)^{12}\right)\\
&=-  4 \pi \lim_{s \to 1} \left(
\sum_{\rho \in H(D)} \calE( \rho, s)  - h   \frac{\Gamma(1/2)
    \Gamma(s-1/2) \zeta(2s-1)}{\Gamma(s) \zeta(2s)}\right) +  12 h
\log(4 \pi)\\
&=   12  h \left( \frac{L'(\chi_D,0)}{L(\chi_D,0)} + \frac{1}{2}
  \log|D| 
\right).
 \end{align*}
Since $m_\rho = 2/w$, the modular form $\Psi_D$ has weight $24 h/w$
and therefore 
\begin{align*}
\frac{3}{\pi} \int\limits_{\Gamma_1\setminus \H} \log\left( |\Psi_D(\tau)| (4 \pi
  y)^{k/2}\right) \frac{dxdy}{ y^2} =& 
  \frac{24 h}{w} 
   \left(2
  \frac{\zeta'(-1)}{\zeta(-1)}+1 \right) -  \frac{24 h}{w} \left( \frac{L'(\chi_D,0)}{L(\chi_D,0)} + \frac{1}{2}
  \log|D| 
\right).
\end{align*}
The claims follows by resorting  the above terms. 
\end{proof}

We remark that the degree formula 
 \eqref{eq:classnumber} is in this case 
the well known analytic class number formula.

Jan Hendrik Bruinier, Mathematisches Institut, Universit\"at Heidelberg, Im Neuenheimer Feld 288, 69120 Heidelberg, Germany. E-mail: bruinier@mathi.uni-heidelberg.de

Ulf K\"uhn, Humboldt-Universit\"at, Institut f\"ur Mathematik, Unter den Linden, 10099 Berlin, Germany. E-mail: kuehn@mathematik.hu-berlin.de

\end{document}